\documentclass[a4paper,12pt]{article}

\usepackage[english]{babel}
\usepackage{amsmath,amsfonts,amssymb,amsthm,bbm,graphics,epsfig,psfrag}
\usepackage{bbm}
\usepackage{paralist}
\usepackage[active]{srcltx}
\usepackage{color,soul}
\setstcolor{red}

\newtheorem{theorem}{Theorem}[section]
\newtheorem{corollary}[theorem]{Corollary}
\newtheorem{lemma}[theorem]{Lemma}
\newtheorem{proposition}[theorem]{Proposition}
\newtheorem{definition}[theorem]{Definition}
\newtheorem{remark}[theorem]{Remark}

\newtheorem{conjecture}[theorem]{Conjecture}

\def\N{\mathbb{N}}

\def\R{\mathbb{R}}
\def\epsilon{\varepsilon}

\let\vp=\varphi

\let\.=\cdot
\let\0=\emptyset

\def\eq#1{{\rm(\ref{eq:#1})}}

\def\tilde{\widetilde}

\newenvironment{formula}[1]{\begin{equation}\label{eq:#1}}{\end{equation}\noindent}
\def\Fi#1{\begin{formula}{#1}}
\def\Ff{\end{formula}\noindent}

\newcommand{\be}{\begin{equation}}
\newcommand{\ee}{\end{equation}}
\newcommand{\baa}{\begin{array}}
\newcommand{\eaa}{\end{array}}
\newcommand{\ba}{\begin{eqnarray}}
\newcommand{\ea}{\end{eqnarray}}

\begin{document}
\date{}
\title{\bf{Transition fronts for the Fisher-KPP equation}}
\author{Fran\c cois Hamel$^{\hbox{\small{ a}}}$ and Luca Rossi$^{\hbox{\small{ b
}}}$\thanks{This work has been carried out in the framework of the Labex
Archim\`ede (ANR-11-LABX-0033) and of the A*MIDEX project (ANR-11-IDEX-0001-02),
funded by the ``Investissements d'Avenir" French Government program managed by
the French National Research Agency (ANR). The research leading to these results
has also received funding from the European Research Council under the European
Union's Seventh Framework Programme (FP/2007-2013) / ERC Grant Agreement
n.321186~- 
ReaDi~- Reaction-Diffusion Equations, Propagation and Modelling and from 
Italian GNAMPA-INdAM. Part of this
work was carried out during visits by F.~Hamel to the Departments of Mathematics
of the University of California, Berkeley, of Stanford University and of the
Universit\`a di Padova, whose hospitality is thankfully acknowledged.}\\
\\
\footnotesize{$^{\hbox{a }}$Aix Marseille Universit\'e, CNRS, Centrale Marseille}\\
\footnotesize{Institut de Math\'ematiques de Marseille, UMR 7373, 13453 Marseille, France}\\
\footnotesize{\& Institut Universitaire de France}\\
\footnotesize{$^{\hbox{b }}$Dipartimento di Matematica P.~e A., Universit\`a di Padova, via Trieste 63, 35121 Padova, Italy}}
\maketitle

\begin{abstract}
This paper is concerned with transition fronts for reaction-diffusion equations
of the Fisher-KPP type. Basic examples of transition fronts connecting the
unstable steady state to the stable one are the standard traveling fronts, but
the class of transition fronts is much larger and the dynamics of the solutions
of such equations is very rich. In the paper, we describe the class of
transition fronts and we study their qualitative dynamical properties. In
particular, we characterize the set of their admissible asymptotic past and
future speeds and their asymptotic profiles and we show that the transition
fronts can only accelerate. We also classify the transition fronts in the class
of measurable superpositions of standard traveling fronts.
\end{abstract}


\section{Introduction and main results}\label{intro}

The purpose of this paper is to describe front-like solutions connecting $0$ 
and $1$ for the equation
\Fi{P=f}
u_t=u_{xx}+f(u),\quad t\in\R,\ x\in\R,
\Ff
when $f$ is of Fisher-KPP (for Kolmogorov-Petrovski-Piskunov) 
type \cite{fi,kpp}. Namely, we will assume that
\be\label{hyp-f}\left\{\baa{l}
f:[0,1]\to\R \text{ is of class $C^2$ and concave},\vspace{3pt}\\
f(0)=f(1)=0\hbox{ and }f(s)>0\hbox{ for all }s\in(0,1).\eaa\right.
\ee
A typical example of such a function $f$ is the quadratic nonlinearity 
$f(s)=s(1-s)$ coming from the logistic equation. Equations of the type~\eq{P=f} 
are among the most studied reaction-diffusion equations and they are known to 
be good models to describe the propagation of fronts connecting the 
unstable steady state $0$ and the stable steady state $1$, see 
e.g.~\cite{f,mu,sk}. The solution $u$ typically stands for the density of a 
species invading an open space and the fronts play a fundamental role in 
the description of the dynamical properties of the solutions. This fact was 
already understood in the  pioneering paper~\cite{kpp}, where the 
authors investigated the existence of solutions of~\eq{P=f} of the form 
$u(t,x)=\vp(x-ct)$ with~$c$ constant and $0=\vp(+\infty)<\vp(\xi)<\vp(-\infty)=1$ 
for all $\xi\in\R$. We will refer to such solutions as {\em standard traveling 
fronts} (connecting~$0$ and~$1$).

For Fisher-KPP functions $f$ satisfying \eqref{hyp-f}, a standard traveling front 
$u(t,x)=\vp(x-ct)$ exists if and only if
$$c\ge2\sqrt{f'(0)}=:c^*$$
and, for each $c\ge c^*$, the function $\vp=\vp_c$ is decreasing and unique up 
to shifts~\cite{aw,f,kpp,u}. Furthermore, these fronts~$\vp_c(x-ct)$ are known 
to be stable with respect to perturbations in some suitable weighted spaces and 
to attract the solutions of the associated Cauchy problem for a large class 
of exponentially decaying initial conditions, see 
e.g.~\cite{bmr,b,ev,hnrr,k,kpp,l,sa,u}. For a given $c\in[c^*,+\infty)$, the 
standard traveling front $\vp_c(x-ct)$ moves with the constant speed~$c$, in 
the sense that its profile $\vp_c(x)$ is translated in time towards right with 
speed~$c$. In particular, its level sets move with speed~$c$. The existence of 
a continuum of speeds is in sharp contrast with the uniqueness of the speed 
for bistable type nonlinearities~$f(u)$, see e.g.~\cite{aw,fm}. The standard 
traveling fronts do not exist in general if the coefficients of the equation depend 
on space or time. When the dependence is periodic, the notion which 
immediately generalizes that of standard traveling fronts is that of 
pulsating traveling fronts, see e.g. \cite{bh,hr2,lz1,lz2,n,nrx,skt,w} for 
the existence and stability properties of such pulsating traveling fronts for time 
and/or space periodic KPP type equations. The notion of pulsating traveling 
front can also be extended in order to deal with time almost-periodic,
almost-automorphic, 
recurrent or uniquely ergodic media, see~\cite{hs,s3,s4,s5,s6,s7}.

However, even for the homogeneous equation~\eq{P=f}, what was 
unexpected until quite recently was the existence of front-like solutions 
connecting $0$ and $1$ which are not standard traveling fronts and whose 
speed changes in time. More precisely, it was proved in~\cite{hn2} that, for any 
real numbers $c_-$ and~$c_+$ such that $2\sqrt{f'(0)}=c^*\le c_-<c_+$, there 
exist solutions~$0<u<1$ satisfying $u(t,-\infty)=1$ and $u(t,+\infty)=0$ 
locally uniformly in $t\in\R$, together with
\be\label{transitionhn}\left\{\baa{ll}
u(t,x)-\vp_{c_-}(x-c_-t)\to0 & \hbox{as }t\to-\infty,\vspace{3pt}\\
u(t,x)-\vp_{c_+}(x-c_+t)\to0 & \hbox{as }t\to+\infty,\eaa\right.\hbox{ uniformly 
in }x\in\R,
\ee
where $\vp_{c_-}(x-c_-t)$ and $\vp_{c_+}(x-c_+t)$ are any two given 
standard traveling fronts with speeds $c_-$ and $c_+$ respectively. 
Roughly speaking, \eqref{transitionhn} implies that the level sets of $u$ 
move approximately with speed $c_-$ for large negative times and with speed 
$c_+$ for large positive times, and this is why one can intuitively say that the 
speed of $u$ changes in time.


\subsubsection*{Transition fronts}

The non-uniformly translating front-like solutions 
satisfying~(\ref{transitionhn}) fall into a notion of fronts connecting 
$0$ and $1$ introduced some years after \cite{hn2}. This new notion,
called transition fronts 
after~\cite{bh1,bh2}, has also been introduced in order to cover the case of
general 
space-time dependent equations. When applied to a one-dimensional 
equation such as~\eq{P=f}, the definition of transition fronts reads as
follows: 

\begin{definition}\label{def1} For problem~\eq{P=f}, a transition front connecting 
$0$ and~$1$ is a time-global classical solution $u:\R\times\R\to[0,1]$ for 
which there exists a function~$X:\R\to\R$ such that
\be\label{gtf}\left\{\baa{ll}
u(t,X(t)+x)\to 1 & \hbox{as }x\to-\infty,\vspace{3pt}\\
u(t,X(t)+x)\to 0 & \hbox{as }x\to+\infty,\eaa\right.\hbox{ uniformly in }t\in\R.
\ee
\end{definition}

The real numbers $X(t)$ therefore reflect the positions of a transition front as 
time runs. The simplest examples of transition fronts connecting $0$ and $1$ 
are the standard traveling fronts $u(t,x)=\vp_c(x-ct)$, for which one can set 
$X(t)=ct$ for all $t\in\R$ (in time or space-periodic media, pulsating traveling 
fronts connecting $0$ and $1$ are also transition fronts and the 
positions~$X(t)$ can be written in the same form). Furthermore,  the solutions 
$u$ constructed in~\cite{hn2} and satisfying~(\ref{transitionhn}) are transition 
fronts connecting $0$ and $1$, with $X(t)=c_-t$ for~$t< 0$ and~$X(t)=c_+t$ 
for~$t\ge0$. This example shows that the above definition of transition fronts 
is meaningful and does not reduce to standard traveling fronts even in the case 
of the homogeneous equation~\eq{P=f}.

Definition \ref{def1} is actually a particular case of a more general definition given
in~\cite{bh1,bh2} in a broader framework. For the one-dimensional 
equation~\eq{P=f}, the transition fronts connec\-ting~$0$ and $1$ correspond to 
the ``wave-like" solutions defined in~\cite{s3,s4} (see also~\cite{m} for a 
different notion involving the continuity with respect to the environment around 
the front position). Roughly speaking, property~(\ref{gtf}) means that the 
diameter of the transition zone between the sets where $u\simeq 1$ and 
$u\simeq 0$ is uniformly bounded in time; the position $X(t)$ is at bounded 
Hausdorff distance from this transition zone. The fundamental property 
in~(\ref{gtf}) is the uniformity of the limits with respect to time~$t\in\R$ (we point 
out that there are solutions satisfying these limits pointwise in $t\in\R$, but not 
uniformly, as shown by Corollary~\ref{cor3} below). Another way of 
writing~(\ref{gtf}) is to say that, given any real numbers $a$ and $b$ such 
that~$0<a\le b<1$, there is a constant~$C=C(u,a,b)\ge0$ such that, for every 
$t\in\R$, 
$$\big\{x\in\R;\ a\le u(t,x)\le b\big\}\subset\big[X(t)-C,X(t)+C].$$
It is easily seen that, for every transition front $u$ connecting $0$ and $1$ and for 
every $C\ge0$,
\be\label{infsup}
0<\inf_{t\in\R,\,x\in[X(t)-C,X(t)+C]}u(t,x)\le\sup_{t\in\R,\,x\in[X(t)-C,X(t)+C]
}u(t,x)<1.
\ee
Indeed, if there were a sequence $(t_n,x_n)_{n\in\N}$ in $\R^2$ such that 
$(x_n-X(t_n))_{n\in\N}$ is bounded and~$u(t_n,x_n)\to0$ as $n\to+\infty$, 
then standard parabolic estimates would imply that the 
functions~$u_n(t,x)=u_n(t+t_n,x+x_n)$ would converge locally uniformly, up 
to extraction of a subsequence, to a solution $0\le u_{\infty}\le1$ of~\eq{P=f} such 
that $u_{\infty}(0,0)=0$, whence~$u_{\infty}=0$ in $\R^2$ from the strong maximum 
principle and the uniqueness of the bounded solutions of the Cauchy problem 
associated to~\eq{P=f}. But~\eqref{gtf} and the boundedness of the 
sequence~$(x_n-X(t_n))_{n\in\N}$ would yield the existence of $M>0$ such that 
$u(t_n,x_n-M)\ge1/2$ for all $n\in\N$, whence~$u_{\infty}(0,-M)\ge1/2$, which is a 
contradiction. Therefore, the first inequality of~(\ref{infsup}) follows, and the last 
one can be proved similarly.

Notice also that, for a given transition front $u$, the family $(X(t))_{t\in\R}$ is not 
uniquely defined since, for any bounded function $\xi:\R\to\R$, the family 
$(X(t)+\xi(t))_{t\in\R}$ satisfies~(\ref{gtf}) if~$(X(t))_{t\in\R}$ does.
Conversely, 
if $(X(t))_{t\in\R}$ and $(\tilde{X}(t))_{t\in\R}$ are associated with a given transition 
front~$u$ connecting~$0$ and $1$ in the sense of~(\ref{gtf}), then it follows 
from~(\ref{gtf}) applied to the family $(X(t))_{t\in\R}$ and from~(\ref{infsup}) applied 
to the family~$(\tilde{X}(t))_{t\in\R}$ that
\be\label{tildeX}
\sup_{t\in\R}\big|X(t)-\tilde{X}(t)\big|<+\infty.
\ee
In other words, we can say that the family of positions $(X(t))_{t\in\R}$ of a given 
transition front connecting~$0$ and $1$ is defined up to an additive bounded
function.

The main goals of the paper are: describing the set of transition fronts for 
equation~\eq{P=f}, providing conditions for a solution $0<u<1$ to be a transition 
front, studying the qualitative properties of the transition fronts, characterizing their 
asymptotic propagation rates for large negative and positive times, in a sense to be 
made more precise below, and finally describing their asymptotic profiles. 


\subsubsection*{Boundedness of the local oscillations of the positions of a transition front}

First of all, we claim that a given transition front cannot move infinitely fast. More 
precisely, even if the positions $(X(t))_{t\in\R}$ of a given transition front are only 
defined up to bounded functions, their local oscillations are uniformly bounded, in 
the sense of the following proposition.

\begin{proposition}\label{pro1} Under the assumptions \eqref{hyp-f}, for any 
transition front of~\eq{P=f} connecting~$0$ and $1$ in the sense of 
Definition~$\ref{def1}$, there holds
\be\label{Xtau}
\forall\,\tau\ge0,\ \ \sup_{(t,s)\in\R^2,\,|t-s|\le\tau}|X(t)-X(s)|<+\infty.
\ee
\end{proposition}

Because this result is of independent interest, we will prove it under weaker 
hypotheses on $f$, in the case of time-space dependent equations, see 
Proposition~\ref{pro1bis} in Section~\ref{sec44} below. 
Applying recursively \eqref{Xtau} with, say, $\tau=1$,
one deduces that 
\be\label{speedbound}
\limsup_{t\to\pm\infty}\left|\frac{X(t)}{t}\right|<+\infty,
\ee
that is, the propagation speeds of any transition front connecting~$0$ and~$1$ are 
asymptotically bounded as $t\to\pm\infty$. This property is in sharp 
contrast with some known acceleration phenomena and infiniteness of the 
asymptotic speed as $t\to+\infty$ for the solutions of the associated Cauchy 
problem when the initial condition~$u_0(x)$ decays to~$0$ as~$x\to+\infty$ more 
slowly than any exponential function: we refer to~\cite{hr1} for the existence of 
such accelerating solutions, which are not transition fronts connecting $0$ and $1$ 
on the time interval~$[0,+\infty)$. These observations somehow mean that the 
transition fronts $u$ connecting $0$ and $1$, as defined by~(\ref{gtf}), cannot 
decay too slowly to $0$ as~$x\to+\infty$, for any given time $t\in\R$.


\subsubsection*{Global mean speed}

As it has become clear from the previous considerations, important notions 
associated with the transition fronts are those of their propagation speeds. We first 
define the notion of possible global mean speed.

\begin{definition} We say that a transition front connecting~$0$ and $1$ 
for~\eq{P=f} has a {\em global mean speed} $\gamma\in\R$ if
\be\label{meanspeed}
\frac{X(t)-X(s)}{t-s}\to\gamma\ \hbox{ as }t-s\to+\infty,
\ee
that is $(X(t+\tau)-X(t))/\tau\to\gamma$ as $\tau\to+\infty$ uniformly in $t\in\R$.
\end{definition}

We point out that this definition slightly differs from the one used in~\cite{bh1,bh2}, 
where absolute values were taken in~(\ref{meanspeed}), that is 
$|X(t)-X(s)|/|t-s|\to|\gamma|$ as $t-s\to+\infty$ (a more general notion of global 
mean speed involving the Hausdorff distance of family of moving hypersurfaces is 
defined in~\cite{bh1,bh2} for more general equations in higher dimensions). But, in 
this one-dimensional situation considered here, it is quite easy to see that the two 
definitions coincide. Notice that, due to Proposition~\ref{pro1}, a given transition 
front connecting $0$ and~$1$ for~\eq{P=f} cannot have an infinite global mean 
speed, that is $\gamma$ in~(\ref{meanspeed}) cannot be equal to~$\pm\infty$. 
Furthermore, if a transition front has a global mean speed $\gamma$, then this 
speed does not depend on the family~$(X(t))_{t\in\R}$, due to the 
property~(\ref{tildeX}). But the global mean speed, if any, does depend on the 
transition front. Actually, any standard traveling front~$\vp_c(x-ct)$ 
sol\-ving~\eq{P=f} has a global mean speed equal to $c$, and the set of all 
admissible global mean speeds among all standard traveling fronts is therefore 
equal to the semi-infinite interval~$[c^*,+\infty)$ with~$c^*=2\sqrt{f'(0)}$. 

Moreover, it turns out that the set of global mean speeds is not increased 
when passing from the set of standard traveling fronts to the larger class of 
transition fronts, and the standard traveling fronts are the only transition fronts with 
non-critical (that is, $\gamma>c^*$) global mean speed~$\gamma$. Indeed, we 
derive the following.

\begin{theorem}\label{thmeanspeed} Under the assumptions \eqref{hyp-f}, the set 
of admissible global mean speeds of transition fronts connecting $0$ and $1$ 
for~\eq{P=f} is equal to the interval $[c^*,+\infty)$. Furthermore, if a transition front 
$u$ connecting $0$ and $1$ has a global mean speed $\gamma>c^*$, then it is a 
standard traveling front of the type $u(t,x)=\varphi_{\gamma}(x-\gamma t)$.
\end{theorem}

In Theorem~\ref{thmeanspeed}, the fact that transition fronts cannot have global 
mean speeds less than~$c^*$ easily follows from standard spreading results of 
Aronson and Weinberger~\cite{aw}, see~(\ref{awspreading}) below. The second 
statement of Theorem~\ref{thmeanspeed}, classifying transition fronts with 
non-minimal global mean speed, is more intricate. It follows from additional results 
of~\cite{hn2} on the classification of solutions of~\eq{P=f} which move faster than 
$c^*$ as~$t\to-\infty$, in a sense to be made more precise later, and from a more 
detailed analysis of the convergence properties of such solutions, see 
Theorem~\ref{thspeeds} below. It is conjectured that the transition fronts with global 
mean speed~$\gamma=c^*$ are also standard traveling 
fronts~$\varphi_{c^*}(x-c^*t)$, but the proof of this conjecture would require further 
arguments.


\subsubsection*{Asymptotic past and future speeds and limiting profiles}

The aforementioned example~(\ref{transitionhn}) shows that, even in a 
homogeneous medium, transition fronts may not admit a global mean speed, since 
the transition fronts~$u$ satisfying~\eqref{transitionhn} are such 
that~$X(t)/t\to c_-$ as $t\to-\infty$ and $X(t)/t\to c_+$ as $t\to+\infty$, with 
$c_+>c_-$. This fact leads us naturally to define the notion of possible asymptotic 
speeds as $t\to-\infty$ and as~$t\to+\infty$.

\begin{definition} We say that a transition front connecting $0$ and~$1$ for the 
equation~\eq{P=f} has an {\em asymptotic past speed} $c_-\in\R$, resp.~an {\em 
asymptotic future speed} $c_+\in\R$, if 
\be\label{cpm2}
\frac{X(t)}{t}\to c_-\hbox{ as }t\to-\infty,\ \hbox{ resp.~}\frac{X(t)}{t}\to c_+\hbox{ as }t\to+\infty.
\ee
\end{definition}

As in the case of the global mean speed, the asymptotic past and future speeds, if 
any, of a given transition front do not depend on the family~$(X(t))_{t\in\R}$. 
Clearly, if a transition front admits a global mean speed $\gamma$ in the sense 
of~(\ref{meanspeed}), then it has asymptotic past and future speeds equal 
to~$\gamma$. Notice also that Proposition~\ref{pro1} implies that a given
transition cannot have past or future speed equal to $\pm\infty$.

The asymptotic past and future speeds, if any, characterize the rate of expansion 
of a given transition front at large negative or positive times. These asymptotic 
speeds might not exist a priori and one could wonder whether these notions of 
speeds as $t\to\pm\infty$ would be sufficient to describe the large time dynamics of 
all transition fronts connecting $0$ and~$1$ for~\eq{P=f}. As a matter of fact, one of 
the main purposes of the present paper will be to characterize completely the set of 
all admissible asymptotic speeds and to show that the asymptotic speeds exist and 
characterize the asymptotic profiles of the fronts, at least for all fronts which are 
non-critical as~$t\to-\infty$.

We first begin with the complete characterization of the set of admissible 
asymptotic past and future mean speeds.

\begin{theorem}\label{thasymptotic}
Under the assumptions \eqref{hyp-f}, transition fronts connecting $0$ and $1$ 
for~\eq{P=f} and having asymptotic past and future speeds $c_{\pm}$ as 
$t\to\pm\infty$ exist if and only if
$$c^*\le c_-\le c_+<+\infty.$$
\end{theorem}

In this theorem, the sufficiency condition is clear. Indeed, firstly, for any 
$c\in[c^*,+\infty)$, the standard traveling fronts $\varphi_c(x-ct)$ are transition 
fronts with asymptotic past and future speeds equal to $c$. Secondly, for any 
$c^*\le c_-<c_+<+\infty$, the transition fronts solving~\eqref{transitionhn} have 
asymptotic past and future speeds equal to $c_{\pm}$. However, unlike the 
Liouville-type result stated in Theorem~\ref{thmeanspeed} for the fronts with
non-critical global mean speeds, transition fronts having prescribed asymptotic
past and future speeds $c_-<c_+$ in $[c^*,+\infty)$ are not unique up to shifts in 
time or space: as a matter of fact, there exists an infinite-dimensional manifold of 
solutions satisfying~(\ref{transitionhn}) for any prescribed speeds $c_-<c_+$ in 
$[c^*,+\infty)$, as will be seen later.

The new part in Theorem~\ref{thasymptotic} is the necessity condition, that is 
transition fronts connecting $0$ and~$1$ and having asymptotic past and future 
speeds can only accelerate and never decelerate, in the sense that the asymptotic 
future speed cannot be smaller than the past one. This property can have the 
following heuristic explanation: if one believes that a transition front $u$ is 
a superposition of a family of standard fronts 
$(\varphi_c(x-ct))_c$ (and this is generally the case, in a sense 
that will be explained later on in this section), then this family should not 
contain fronts with $c$ less than the asymptotic past speed $c_-$ of $u$ 
(otherwise, since they are slower, they would be dominant as~$t\to-\infty$).
Therefore, for very negative 
times~$t$, the tail of $u(t,\cdot)$ at~Ê$+\infty$ should not involve the faster decay 
rates~$\lambda_c$, with $\lambda_c>\lambda_{c_-}$, of the traveling fronts with 
smaller speeds $c<c_-$. In other words, for very negative times~$t$, $u(t,\cdot)$ is 
expected to decay at~$+\infty$ with a rate not larger than $\lambda_{c_-}$. But 
since such tails give rise in future times to solutions which move not slower than 
$c_-$, by~\cite{u}, $u$ is then expected to move faster than~$c_-$ as 
$t\to+\infty$.

Actually, the necessity condition of Theorem~\ref{thasymptotic} is a consequence 
of the following more general and more precise result.

\begin{theorem}\label{thliminf}
Under the assumptions \eqref{hyp-f}, for any transition front $u$ connecting $0$ 
and $1$ for~\eq{P=f}, there holds
\be\label{inequalities}
c^*\le\liminf_{t\to-\infty}\frac{X(t)}{t}\le\liminf_{t\to+\infty}\frac{X(t)}{t}\le
\limsup_{t\to+\infty}\frac{X(t)}{t}<+\infty.
\ee
Furthermore, if $c^*<\liminf_{t\to-\infty}X(t)/t$, then $u$ has asymptotic past and 
future speeds, that is the limits $c_{\pm}=\lim_{t\to\pm\infty}X(t)/t$ exist in $\R$, 
with $c_-\le c_+$, and there exists a bounded function $\xi:\R\to\R$ such that
\be\label{cvfronts0}
u(t,X(t)+\xi(t)+\cdot)\to\varphi_{c_{\pm}}\ \hbox{ in }C^2(\R)\ \hbox{ as }t\to\pm\infty.
\ee
\end{theorem}

Some comments are in order. The first inequality of~(\ref{inequalities}) will follow 
immediately from some spreading estimates of~\cite{aw}, see~\eqref{awspreading} 
below. The last inequality is a direct consequence of Proposition~\ref{pro1}. The 
difficult inequality is the second one. Actually, if~$c^*=\liminf_{t\to-\infty}X(t)/t$, then 
the inequality $\liminf_{t\to-\infty}X(t)/t=c^*\le\liminf_{t\to+\infty}X(t)/t$ also follows 
from~\cite{aw}. Therefore, the real interest of Theorem~\ref{thliminf} lies in the 
second assertion. Among other things, its proof uses the decomposition of 
supercritical transition fronts as superposition of standard traveling fronts. This 
decomposition, given in~\cite{hn2}, will be explained at the end of this section and it 
is conjectured to hold for {\it every} solution $0<u<1$ of~\eq{P=f}. Therefore, we 
conjecture that the following property holds, even in the critical case 
where~$\liminf_{t\to-\infty}X(t)/t=c^*$.

\begin{conjecture}
Under the assumptions \eqref{hyp-f}, any transition front $u$ connecting $0$ and 
$1$ for~\eq{P=f} has asymptotic past and future speeds $c_{\pm}$ such that 
$c^*\le c_-\le c_+<+\infty$ and the convergence~\eqref{cvfronts0} to the limiting 
profiles $\varphi_{c_{\pm}}$ holds.
\end{conjecture}

Theorem~\ref{thliminf} shows the existence of the asymptotic past and future 
speeds for transition fronts $u$ such that $\liminf_{t\to-\infty}X(t)/t$ is not critical 
(that is, larger than $c^*$). Thus, for such fronts, some complex dynamics where 
the set of limiting values of $X(t)/t$ as $t\to-\infty$ and~$t\to+\infty$ would not be 
reduced to a singleton are impossible. As a matter of fact, such infinitely many 
oscillations are possible in general for the solutions of the Cauchy problem 
associated with~\eq{P=f}, as shown in~\cite{hn1,y}. This occurs when the initial 
condition~$u_0(x)$ has itself a complex behavior as~$x\to+\infty$ with infinitely 
many oscillations on some larger and larger intervals between two exponentially 
decaying functions (see also~\cite{ggn} for the existence of asymptotic oscillations 
between two different speeds for the solutions of $x$-dependent Cauchy problems 
when the initial condition~$u_0(x)$ is a Heaviside-type function and $f(x,u)$ is 
slowly and non-periodically oscillating as $x\to+\infty$). Here, we consider
time-global solutions defined for all $t\in\R$ and, roughly speaking, 
Theorem~\ref{thliminf} means that, at least when the transition front is not critical as 
$t\to-\infty$, the possible oscillations of the tails as $x\to+\infty$ are somehow 
annihilated by the fact that the solution is ancient. 


\subsubsection*{Characterization of transition fronts and their asymptotic past and future speeds}

Transition fronts connecting $0$ and $1$ for~\eq{P=f} are particular time-global 
solutions $0<u<1$ satisfying~(\ref{gtf}) for some family $(X(t))_{t\in\R}$. One can 
now wonder what condition would guarantee a given time-global solution $u$ 
of~\eq{P=f} to be a transition front. First of all, the spreading results of~\cite{aw} 
imply that any solution~$0<u(t,x)<1$ of~\eq{P=f} is such that
\be\label{awspreading}
\forall\,c\in[0,c^*),\ \ \max_{[-c|t|,c|t|]}u(t,\cdot)\to0\hbox{ as }t\to-\infty.
\ee
The next result provides a characterization of the transition fronts connecting $0$ 
and $1$ in the class of solutions $0<u<1$ of~\eq{P=f} which satisfy a condition 
slightly stronger than~(\ref{awspreading}).

\begin{theorem}\label{thm:decay}
Assume that $f$ satisfies \eqref{hyp-f} and let $0<u<1$ be a solution of~\eq{P=f} 
such that
\be\label{hyphn}
\exists\,c>c^*,\ \ \max_{[-c|t|,c|t|]}u(t,\cdot)\to0\ \hbox{ as }t\to-\infty.
\ee
Then the limit
\Fi{decay}
\lambda:=-\lim_{x\to+\infty}\frac{\ln u(0,x)}{x}
\Ff
exists, satisfies $\lambda\in\big[0,\sqrt{f'(0)}\big)$ and $u$ is a
transition front connecting 
$0$ and $1$ if and only if~$\lambda>0$. Furthermore, if $\lambda>0$, the transition 
front $u$ admits asymptotic past and future speeds $c_-$ and $c_+$ given by
\be\label{asympspeeds}
c^*<c_-=\sup\Big\{c\ge0,\ \lim_{t\to-\infty}\max_{[-c|t|,c|t|]}u(t,\cdot)=0\Big\}\leq 
c_+=\lambda+\frac{f'(0)}\lambda,
\ee
and~\eqref{cvfronts0} holds for some bounded function $\xi:\R\to\R$.
\end{theorem}

The time $t=0$ in condition~\eq{decay} can obviously be replaced by an arbitrary 
time~$t_0\in\R$. The fact that it is sufficient to know at least some informations on 
the profile of $u$ at a given time in order to infer what its asymptotic future speed 
$c_+$ will be is not surprising, due to the parabolic maximum principle. A
formula for the past speed $c_-$ which would be as simple as the one for $c_+$
is far to be trivial: our characterization of $c_-$ requires a knowledge of
the profile of $u$ as $t\to-\infty$, but it seems hopeless to be able to derive
the past speed from the profile at a given time.

A consequence of Theorem~\ref{thm:decay} is that, under the assumption~\eqref{hyphn}, if $0<u<1$ solves~\eq{P=f} but is not a transition front connecting $0$ and $1$, then, for every $\epsilon>0$, there is a constant $C_{\epsilon}>0$ such that $u(0,x)\ge C_{\epsilon}\,e^{-\epsilon x}$ for all $x\ge0$. In other words, for any solution~$u$ which is not a transition front and for any time, say $t=0$, then~$u(0,x)$ either decays to~$0$ as~$x\to+\infty$ more slowly than any exponentially decaying function, or is away from $0$ at~$+\infty$ (in the latter case, it follows from the proof of Theorem~\ref{thm:decay} that $\lim_{x\to+\infty}u(0,x)$ exists in~$(0,1]$).

Hypothesis \eqref{hyphn} is required in the proof of Theorem~\ref{thm:decay} in 
order to apply the result of classification of such solutions $0<u<1$ of~\eq{P=f}, 
see~\cite{hn2} and the end of this section. If, as conjectured, the characterization 
holds without this additional condition~\eqref{hyphn}, then we can reasonably state 
the following conjecture.

\begin{conjecture}
Under the assumptions \eqref{hyp-f}, for any solution $0<u<1$ of~\eq{P=f}, the 
limit~$\lambda=-\lim_{x\to+\infty}\ln u(0,x)/x$ exists in $[0,\sqrt{f'(0)}]$ and $u$ is a 
transition front connecting~$0$ and $1$ if and only if $\lambda>0$. Furthermore, if 
$\lambda>0$, then the transition front $u$ has asymptotic past and future speeds 
$c_-$ and $c_+$, which satisfy~\eqref{asympspeeds} with $c^*\le c_-$, 
and~\eqref{cvfronts0} holds.
\end{conjecture}

Finally, we point out that Theorem~\ref{thm:decay} will easily imply the second part 
of the conclusion of Theorem~\ref{thliminf} since, for a given transition front 
connecting $0$ and $1$ such that~$\liminf_{t\to-\infty}X(t)/t>c^*$, the 
condition~\eqref{hyphn} is automatically fulfilled.


\subsubsection*{Transition fronts as superposition of standard traveling fronts}
\label{sec:hn2}

An important tool in the proof of the above theorems is the construction of time-
global solutions $0<u<1$ of~\eq{P=f} as superposition of standard traveling fronts. 
We recall that $\varphi_c(x-ct)$ denotes a standard traveling front for~\eq{P=f} 
such that~$\varphi_c(-\infty)=1>\varphi_c>\varphi_c(+\infty)=0$. Such front exists 
if and only if $c\ge c^*=2\sqrt{f'(0)}$. Furthermore, for each $c\ge c^*$, the 
function~$\varphi_c$ is decreasing, unique up to shift, and one can assume without 
loss of generality that it satisfies
\be\label{asymvarphi}\left.\begin{cases}
\varphi_c(\xi)\sim e^{-\lambda_c\xi}\ \ \hbox{ if }c>c^*,\\
\varphi_{c^*}(\xi)\sim \xi\,e^{-\lambda_{c^*}\xi},
\end{cases}\right.\hbox{ as }\xi\to+\infty,
\ee
where
\be\label{lambdac}
\lambda_c=\frac{c-\sqrt{c^2-4f'(0)}}{2}\quad\hbox{for }c\ge c^*.
\ee
Notice in particular that $\lambda_{c^*}=c^*/2=\sqrt{f'(0)}$. With the 
normalization~(\ref{asymvarphi}), it is known that $\varphi_c(\xi)\le 
e^{-\lambda_c\xi}$ for all $\xi\in\R$ if $c>c^*$. Lastly, let $\theta:\R\to(0,1)$ be the 
unique solution of
$$\theta'(t)=f(\theta(t)),\ t\in\R,$$
such that $\theta(t)\sim e^{f'(0)t}$ as $t\to-\infty$.

Let us now introduce some notations which are used for the construction of other 
solutions of~\eq{P=f} given as measurable interactions of some shifts of the space-
independent solution~$\theta(t)$ and the right- or left-moving standard traveling 
fronts~$\varphi_c(\pm x-ct)$. More precisely, let 
first~$\Psi:[-\lambda_{c^*},\lambda_{c^*}]\to 
X:=(-\infty,-c^*]\cup[c^*,+\infty)\cup\{\infty\}$ be the bijection defined by
$$\left\{\baa{rcll}
\Psi(\lambda) & = & \displaystyle\lambda+\frac{f'(0)}{\lambda} & \hbox{if }
\lambda\in[-\lambda_{c^*},\lambda_{c^*}]\backslash\{0\},\vspace{3pt}\\
\Psi(0) & = & \infty,&\eaa\right.$$
and let us endow $X$ with the topology induced by the image by $\Psi$ of the Borel 
topo\-logy of~$[-\lambda_{c^*},\lambda_{c^*}]$. In other words, a subset $O$ of 
$X$ is open if $\Psi^{-1}(O)$ is open relatively in $[-\lambda_{c^*},\lambda_{c^*}]$, 
that is the intersection of $[-\lambda_{c^*},\lambda_{c^*}]$ with an open set of $\R$. 
The bijection $\Psi^{-1}:X\to[-\lambda_{c^*},\lambda_{c^*}]$ denotes the reciprocal 
of the function $\Psi$ and it is given by~$\Psi^{-1}(c)=\lambda_c:=
(|c|-\sqrt{c^2-4f'(0)})\times({\rm{sgn}}\,c)/2$ if $c\in X\backslash\{\infty\}$ and 
$\Psi^{-1}(\infty)=0$, where~${\rm{sgn}}\,c=c/|c|$ denotes the sign of $c$. Notice 
that the notation for $\lambda_c$ is coherent with~(\ref{lambdac}) when 
$c\in[c^*,+\infty)$. Lastly, let $\mathcal{M}$ be the set of all nonnegative Borel 
measures~$\mu$ on~$X$ such that~$0<\mu(X)<+\infty$.

It follows from Theorem~1.2 of~\cite{hn2} and formula~(30) of~\cite{hn2} that there 
is a one-to-one map
$$\mu\mapsto u_{\mu}$$
from $\mathcal{M}$ to the set of solutions $0<u<1$ of~\eq{P=f}. Furthermore, for 
each~$\mu\in\mathcal{M}$, denoting $M=\mu\big(X\backslash\{-c^*,c^*\}\big)$, the 
solution $u_{\mu}$ constructed in~\cite{hn2} satisfies
\be\label{umu}\baa{l}
\max\Big(\varphi_{c^*}\big(x-c^*t-c^*\ln\mu(c^*)\big),\,\varphi_{c^*}\big(\!-x-c^*t-c^*\ln\mu(-c^*)\big),\vspace{3pt}\\
\qquad\displaystyle\ \ M^{-1}\!\!\int_{\R\backslash[-c^*,c^*]}\!\!\varphi_{|c|}\big(({\rm{sgn}}\,c)x\!-\!|c|t\!-\!|c|\ln M\big)\,d\mu(c)\,+\,M^{-1}\theta(t\!+\!\ln M)\,\mu(\infty)\Big)\vspace{3pt}\\
\le\,u_{\mu}(t,x)\,\le\,\varphi_{c^*}\big(x-c^*t-c^*\ln\mu(c^*)\big)+\varphi_{c^*}\big(\!-x-c^*t-c^*\ln\mu(-c^*)\big)\vspace{3pt}\\
\qquad\qquad\qquad\displaystyle+M^{-1}\!\!\int_{\R\backslash[-c^*,c^*]}\!\!e^{-\lambda_{|c|}(({\rm{sgn}}\,c)x-|c|t-|c|\ln M)}d\mu(c)\,+\,M^{-1}e^{f'(0)(t+\ln M)}\mu(\infty)\eaa
\ee
for all $(t,x)\in\R^2$, under the convention that 
$\mu(\pm c^*)=\mu\big(\{\pm c^*\}\big)$, that $\mu(\infty)=\mu\big(\{\infty\}\big)$ and 
that the terms involving $\ln 0$ are not present. The estimate~(\ref{umu}) reflects 
different types of contributions weighted by $\mu$: critical standard traveling fronts, 
supercritical standard traveling fronts and spatially homogeneous solutions. More 
precisely, from Section~3 of~\cite{hn2}, each~$u_{\mu}$ is constructed as the 
monotone $C^{1,2}_{loc}(\R\times\R)$ limit as $n\to+\infty$ of the functions 
$u^n_{\mu}$ defined in $[-n,+\infty)\times\R$ as the solutions of the Cauchy
problems
\be\label{defumun}\left\{\baa{rcl}
(u_{\mu}^n)_t & \!\!=\!\! & (u_{\mu}^n)_{xx}+f(u^n_{\mu}),\ \ t>-n,\ x\in\R,\vspace{3pt}\\
u_{\mu}^n(-n,x) & \!\!=\!\! & \max\Big(\varphi_{c^*}\big(x+c^*n-c^*\ln\mu(c^*)\big),\,\varphi_{c^*}\big(\!-x+c^*n-c^*\ln\mu(-c^*)\big),\vspace{3pt}\\
& \!\!\!\! & \qquad\displaystyle\ \ M^{-1}\!\!\int_{\R\backslash[-c^*,c^*]}\!\!\varphi_{|c|}\big(({\rm{sgn}}\,c)x+|c|n-|c|\ln M\big)\,d\mu(c)\vspace{3pt}\\
& \!\!\!\! & \qquad\ \ +M^{-1}\theta(-n+\ln M)\,\mu(\infty)\Big),\ \ x\in\R.\eaa\right.
 \ee
In particular, for each $\mu\in\mathcal{M}$, one has $u^n_{\mu}(t,x)\to u_{\mu}(t,x)$ 
as $n\to+\infty$ locally uniformly in~$\R^2$. Throughout the paper, the notation 
$u_\mu$ will refer to such solutions, which all fulfill~$0<u_{\mu}<1$ in $\R^2$. A 
straightforward consequence of the maximum principle is that the 
solutions~$u_{\mu}$ are nonincreasing (resp.~nondecreasing) with respect to $x$ 
if $\mu\big((-\infty,-c^*]\big)=0$ (resp.~if $\mu\big([c^*,+\infty)\big)=0$), that is if 
$u_{\mu}$ is a measurable interaction of right-moving (resp. left-moving) spatially 
decreasing (resp. increasing) traveling fronts $\varphi_c(x-ct)$ (resp. 
$\varphi_c(-x-ct)$) and the spatially uniform function~$\theta$. Before going further 
on, let us just mention the simplest example of such solutions $u_{\mu}$. Namely, if 
$\mu=m\,\delta_c$ with $m>0$ and $c\in(-\infty,-c^*]\cup[c^*,+\infty)$, 
where~$\delta_c$ denotes the Dirac measure at $\{c\}$, then $u_{\mu}^n(-n,x)=
\varphi_{|c|}\big(({\rm{sgn}}\,c)x+|c|n-|c|\ln m\big)$ for all~$n\in\N$ and~$x\in\R$, whence~$u_{\mu}(t,x)=\varphi_{|c|}\big(({\rm{sgn}}\,c)x-|c|t-|c|\ln m\big)$ for all 
$(t,x)\in\R^2$. For more general measures $\mu$, we refer to~\cite{hn2} for 
additional monotonicity properties, continuity with respect to~$\mu$ and asymptotic 
properties of~$u_{\mu}$ as~$t\to\pm\infty$.

Lastly, we point out that these solutions $u_{\mu}$ almost describe the set of all 
solutions of~\eq{P=f}. Indeed, reminding the general property~\eqref{awspreading}, 
the following almost-uniqueness result was proved in~\cite{hn2}: if a 
solution~$0<u(t,x)<1$ of~\eq{P=f} satisfies~\eqref{hyphn}, that is
$$\exists\,c>c^*,\ \ \max_{[-c|t|,c|t|]}u(t,\cdot)\to0\ \hbox{ as }t\to-\infty,$$
then there is a measure~$\mu\in\mathcal{M}$ such that $u=u_{\mu}$ and the 
support of $\mu$ does not intersect the interval $(-c,c)$. For the one-dimensional 
equation~\eq{P=f}, it is conjectured that this classification result holds in full 
generality, that is, without assuming~(\ref{hyphn}).

This almost-uniqueness result motived us to investigate more deeply the properties 
of the solutions of~\eq{P=f} of the type $u_\mu$, starting from the following very 
natural question: under which condition on $\mu$ is $u_\mu$ a transition front? The 
answer is given by the following statement.

\begin{theorem}\label{thsupport}
Assume that $f$ satisfies \eqref{hyp-f} and let $u_{\mu}$ be the solution 
of~$\eq{P=f}$ associated with a measure~$\mu\in\mathcal{M}$. Then $u_{\mu}$ is 
a transition front connecting $0$ and $1$ if and only if the support of~$\mu$ is 
bounded and is included in $[c^*,+\infty)$. In such a case, $u_{\mu}$ is decreasing 
with respect to the $x$ variable.
\end{theorem}

Theorem~\ref{thsupport} has the following intuitive explanations: since right-moving 
traveling fronts $\varphi_c(x-ct)$ connect $1$ at $-\infty$ and $0$ at $+\infty$, while 
the left-moving ones $\varphi_c(-x-ct)$ connect $0$ at $-\infty$ and $1$ at 
$+\infty$, it is natural that transition fronts connecting $1$ (at~$-\infty$) and $0$ (at 
$+\infty$) only involve right-moving traveling fronts $\vp_c(x-ct)$. Furthermore, the 
faster the standard traveling fronts~$\varphi_c(x-ct)$, the flatter, since the 
quantities $|c\varphi_c'(x-ct)|=\partial_t\varphi_c(x-ct)=
\partial_{xx}\varphi_c(x-ct)+f(\varphi_c(x-ct))$ remain bounded as $c\to+\infty$, by 
standard parabolic estimates. Therefore, for a solution~$u_{\mu}$ of~\eq{P=f} to 
have a uniformly bounded transition region between values close to $1$ and values 
close to $0$, it is natural that the set of speeds of the standard traveling fronts 
involved in the measure~$\mu$ be bounded. 

Based on Theorem~\ref{thsupport}, another important example of transition fronts 
connecting $0$ and~$1$ of the type $u_{\mu}$ is the following one: namely, 
if~$c_-<c_+\in[c^*,+\infty)$ are given, if $\tilde{\mu}\in\mathcal{M}$ is given with 
$\tilde{\mu}\big(X\backslash(c_-,c_+)\big)=0$, then for any positive real numbers 
$m_{\pm}$, the solution $u_{\mu}$ with~$\mu=m_-\delta_{c_-}+m_+\delta_{c_+}+
\tilde{\mu}$ satisfies
$$\sup_{x\in\R}\Big|u_{\mu}(t,x)-\varphi_{c_{\pm}}\big(x-c_{\pm}t-(c_{\pm}-
\lambda_{c_{\pm}}^{-1})\ln(m_-+m_++\tilde{\mu}(X))-
\lambda_{c_{\pm}}^{-1}\ln m_{\pm}\big)\Big|\to0$$
as $t\to\pm\infty$, as follows easily from~\eqref{umu} and Theorem~\ref{thspeeds} 
below; in particular, there are some positive real numbers $m_{\pm}$ such that 
$u_{\mu}$ satisfies~\eqref{transitionhn}. Thus, since two different such measures 
$\tilde{\mu}$ give rise to two different solutions $u_{\mu}$, which are transition 
fronts connecting~$0$ and $1$ by Theorem~\ref{thsupport}, one infers that, for any 
two given speeds $c_-<c_+$ in $[c^*,+\infty)$, there is an infinite-dimensional 
manifold of transition fronts of~\eq{P=f} satisfying~\eqref{transitionhn}, thus having 
asymptotic past and future speeds $c_\pm$.

An immediate byproduct of Theorem~\ref{thsupport}, which is interesting by itself but 
was a priori not obvious, is that the spatial monotonicity together with the existence 
of the limits $0$ and $1$ as~$x\to\pm\infty$ pointwise in $t\in\R$ do not guarantee 
the transition front property~(\ref{gtf}). Indeed, the following corollary holds.

\begin{corollary}\label{cor3}
Under the assumption \eqref{hyp-f}, equation~\eq{P=f} admits some 
solutions~$0<u<1$ such that, for every $t\in\R$, $u(t,\cdot)$ is decreasing with 
$u(t,-\infty)=1$ and $u(t,+\infty)=0$, but which are not transition fronts connecting 
$0$ and $1$.
\end{corollary}

\noindent{\bf{Proof.}} It suffices to take $u=u_{\mu}$, where $\mu$ is any 
measure in $\mathcal{M}$ whose support is included in~$[c^*,+\infty)$ but is not 
compact. Indeed, on the one hand, for such a $\mu$, we know 
from~(\ref{defumun}) and the maximum principle that $u$ is nonincreasing with 
respect to~$x$. Furthermore, the inequa\-lities~(\ref{umu}), together 
with~$\mu\big((-\infty,-c^*]\cup\{\infty\}\big)=0$ and Lebesgue's dominated 
convergence theorem, yield~$u(t,-\infty)=1$ and $u(t,+\infty)=0$ for all $t\in\R$ 
(for the limit $u(t,+\infty)=0$, one uses the upper bound in~(\ref{umu}) and 
the fact that~$0<\lambda_cc=\lambda_c^2+f'(0)\le2f'(0)$ for all~$c\in[c^*,+\infty)$). 
In particular, for any fixed time, the space derivative 
$v=(u_{\mu})_x$ - which is a nonpositive solution of a linear parabolic equation - 
cannot be identically equal to~$0$ and is thus negative in $\R^2$ by the strong 
maximum principle. On the other hand, since~$\mu$ is not compactly 
supported,~$u_{\mu}$ is not a transition front, by 
Theorem~\ref{thsupport}.\hfill$\Box$

\begin{remark}{\rm It actually follows from Theorem~$1.4$ of~$\cite{hn2}$ and from 
part~{\rm{(ii)}} of Theorem~$1.1$ of Zlato{\v{s}}~$\cite{z}$ $($after some changes 
of notations$)$ that, if the measure $\mu\in\mathcal{M}$ is compactly supported 
in~$(c^*,+\infty)$, then~$u_{\mu}$ is a transition front connecting $0$ and $1$ 
for~\eq{P=f}. The paper~$\cite{z}$ is actually concerned with more general 
heterogeneous KPP equations. In the case of the homogeneous equation~\eq{P=f}, 
the proof of the sufficiency condition of Theorem~$\ref{thsupport}$ above is 
actually alternate to that of~$\cite{z}$ and it covers the case where the leftmost 
point of the support of $\mu$ is equal to $c^*$. The proof of 
Theorem~$\ref{thsupport}$, given in Section~$\ref{sec2}$ below, is based on some 
explicit comparisons with some standard traveling fronts (see in particular 
Proposition~$\ref{pro3}$ below), which have their own interest. More precisely, for 
the proof of the sufficiency condition of Theorem~$\ref{thsupport}$, we will show,  
thanks to an intersection number argument, that, if a measure $\mu\in\mathcal{M}$ 
has a bounded support included in~$[c^*,+\infty)$, then the solution $u_{\mu}$ is 
bounded from above (resp. from below) by a suitably shifted standard traveling 
front on the right (resp. on the left) of any point, at any time. On the other hand, the 
proof of the necessity condition of Theorem~\ref{thsupport} is based directly on the 
estimates~(\ref{umu}) and on Proposition~\ref{pro1}.}
\end{remark}

We next derive several properties of the solution $u_\mu$ of~\eq{P=f} in the case it 
is a transition front connecting $0$ and $1$, that is when the measure $\mu$  is 
compactly supported in $[c^*,+\infty)$.

\begin{theorem}\label{thspeeds}
Assume that $f$ satisfies \eqref{hyp-f} and let $u_{\mu}$ be the solution 
of~$\eq{P=f}$ associated with a measure~$\mu\in\mathcal{M}$ whose support is a 
compact subset of $[c^*,+\infty)$. Calling $c_-$ and~$c_+$ the leftmost and 
rightmost points of the support of $\mu$, the following properties hold:
\begin{enumerate}
\item[{\rm{(i)}}] The front $u_{\mu}$ has an asymptotic past speed equal to $c_-$ 
and an asymptotic future speed equal to $c_+$.
\item[{\rm{(ii)}}] The positions $(X(t))_{t\in\R}$ satisfy
$$\limsup_{t\to\pm\infty}|X(t)-c_\pm t|<+\infty\ \text{ if }\mu(c_\pm)>0\ 
\hbox{ and }\ \lim_{t\to\pm\infty}(X(t)-c_\pm t)=-\infty\ \text{ if }\mu(c_\pm)=0.$$
\item[{\rm{(iii)}}]  If $c_->c^*$, then there is a bounded function $\xi:\R\to\R$ such 
that
\be\label{convc+-}
u_{\mu}(t,X(t)+\xi(t)+\cdot)\to\varphi_{c_{\pm}}\ \hbox{ in }C^2(\R)\ \hbox{ as }
t\to\pm\infty.
\ee
\item[{\rm{(iv)}}] If $c_-=c^*$ and~$\mu(c^*)>0$, then 
\be\label{convc*}
u_\mu(t,c^*t+c^*\ln\mu(c^*)+\cdot)\to\varphi_{c^*}\ \hbox{ in }C^2(\R)\ \hbox{ as }
t\to-\infty.
\ee
\end{enumerate}
\end{theorem}

Theorems~\ref{thsupport} and \ref{thspeeds} are the keystones of our paper and 
they are actually the first results we will prove in the sequel. The other theorems 
can be viewed as their corollaries. We chose to present Theorems~\ref{thsupport} 
and~\ref{thspeeds} only at the end of the introduction because their contribution to 
the understanding of general propagation properties for the equation~\eq{P=f} is 
less evident than that of their consequences stated before. Indeed, even though the 
aforementioned quasi-uniqueness result of~\cite{hn2} almost completely 
characterize the class of solutions~$0<u<1$ of~\eq{P=f} in terms of the functions 
$u_\mu$, recovering $\mu$ from $u$ - that is describing the inverse of the 
mapping~$\mu\mapsto u_\mu$ - is a difficult task. Nevertheless, if $u$ 
satisfies~\eqref{hyphn}, Theorem~\ref{thm:decay} above provides a simple 
criterion to to know whether $u$ is a transition front or not and to identify the least 
upper bound of the support of the associated measure $\mu$.

An important point in Theorem~\ref{thspeeds} is the existence and 
characterization of the asymptotic past and future speeds when the support of 
$\mu$ is compactly included in $[c^*,+\infty)$, as well as the characterization of the 
boundedness of the asymptotic shifts $X(t)-c_{\pm}t$ as $t\to\pm\infty$. In 
particular, an immediate consequence of statement~(ii) of Theorem \ref{thspeeds} 
is that, even when the limits $\lim_{t\to\pm\infty}X(t)/t=c_{\pm}$ exist, the distance 
between $X(t)$ and $c_-t$ (resp. $c_+t$) may or may not bounded as $t\to-\infty$ 
(resp. as $t\to+\infty$). Actually, the four possible cases can occur, as it is seen by 
applying Theorem~\ref{thspeeds} with the following measures:
$$\mu_{(\alpha,\beta)}=\alpha\,\delta_{c_-}+\mathbbm{1}_{(c_-,c_+)}+
\beta\,\delta_{c_+},\quad\alpha,\beta\in\{0,1\},$$
where $\mathbbm{1}_{(c_-,c_+)}$ denotes the measure whose density is given by 
the characteristic function of~$(c_-,c_+)$.

In parts~(iii) an~(iv) regarding the convergence of $u_{\mu}$ to the profiles 
$\varphi_{c_{\pm}}$ as~$t\to\pm\infty$, the case $c_-=c^*$ and $\mu(c^*)=0$ is 
not treated (and neither is the case $\mu(c^*)>0$ for the convergence 
as~$t\to+\infty$), the difficulty coming from the fact that the classification of the 
solutions of~\eq{P=f} is still an open question when one does not 
assume~(\ref{hyphn}). If, as conjectured, all solutions of~\eq{P=f} were known to 
be of the type $u_{\mu}$, then~(\ref{convc+-}) would hold without any restriction 
on~$c_-$. However, even the case $c_->c^*$, the proof of the convergence 
formula~(\ref{convc+-}) to the asymptotic profiles $\varphi_{c_{\pm}}$ as 
$t\to\pm\infty$ is far from obvious, due to the fact that in general the quantities 
$X(t)-c_{\pm}t$ are not bounded as $t\to\pm\infty$ and the positions $X(t)$ are not 
those of the traveling fronts $\varphi_{c_{\pm}}(x-c_{\pm}t)$ up to bounded 
functions. The proof of~\eqref{convc+-} is obtained after some passages to the limit 
and the decomposition of the measure $\mu$ into some restricted measures. The 
proof also uses some intersection number arguments.

\begin{remark} {\rm All the above results have a specular counterpart for 
left-moving transition fronts connecting $0$ at $-\infty$ and $1$ at $+\infty$, as it is 
immediately seen by performing the reflection $x\mapsto-x$.}
\end{remark}

\noindent{\bf{Outline of the paper.}} In Section~\ref{sec2}, we will prove 
Theorem~\ref{thsupport} on the characterization of transition fronts given as 
measurable superpositions of standard traveling fronts of the 
type~$\varphi_c(\pm x-ct)$. Section~\ref{sec3} is devoted to the proof of 
Theorem~\ref{thspeeds} on the existence and characterization of the asymptotic 
past and future speeds and the asymptotic profiles of the transition fronts of the 
type $u_{\mu}$. Lastly, we will carry out the proofs of the remaining results, that is 
Theorems~\ref{thmeanspeed},~\ref{thliminf},~\ref{thm:decay} and 
Proposition~\ref{pro1} (in a more general framework), in Section~\ref{sec4}. We 
recall that Theorem~\ref{thasymptotic} essentially follows from 
Theorem~\ref{thliminf}, as explained after the statement of 
Theorem~\ref{thasymptotic}.


\section{Proof of Theorem~\ref{thsupport}}\label{sec2}

We first prove the necessity condition in Section~\ref{sec21}. Section~\ref{sec22} is devoted to the proof of the sufficiency condition. It is based on the auxiliary Proposition~\ref{pro3} below, which is proved in Section~\ref{sec23}.


\subsection{Proof of the necessity condition}\label{sec21}

Here and in the sequel, we use the notation taken from \cite{hn2} and recalled in Section \ref{intro}. Let $\mu\in\mathcal{M}$ be given and assume that~$u_{\mu}$ is a transition front connecting $0$ and $1$ for the equation~\eq{P=f}, that is~$0<u_{\mu}<1$ in~$\R^2$ and there is a family of real numbers $(X(t))_{t\in\R}$ such that~(\ref{gtf}) holds. We will show that the support of $\mu$ is bounded and that it is included in $[c^*,+\infty)$, the latter meaning that~$\mu\big((-\infty,-c^*]\cup\{\infty\}\big)=0$. First of all, if $\mu(-c^*)>0$, then~(\ref{umu}) yields
$$u_{\mu}(0,x)\ge\varphi_{c^*}(-x-c^*\ln\mu(-c^*))\to1\ \hbox{ as }x\to+\infty,$$
which immediately contradicts~(\ref{gtf}). On the other hand, if $\mu\big((-\infty,-c^*)\cup\{\infty\}\big)>0$, then~$M=\mu\big(X\backslash\{-c^*,c^*\}\big)>0$ and it follows from~(\ref{umu}) that
$$u_{\mu}(0,x)\ge M^{-1}\int_{(-\infty,-c^*)}\varphi_{|c|}\big(-x-|c|\ln M\big)\,d\mu(c)\,+\,M^{-1}\theta(\ln M)\,\mu(\infty)$$
for all $x\in\R$. If $\mu(\infty)>0$, then 
$\inf_{\R}u_{\mu}(0,\cdot)\ge M^{-1}\theta(\ln M)\,\mu(\infty)>0$, which
contradicts~(\ref{gtf}). If $\mu\big((-\infty,-c^*)\big)>0$, then applying Lebesgue's dominated convergence theorem, one
deduces that $\liminf_{x\to+\infty}u_{\mu}(0,x)\ge
M^{-1}\mu\big((-\infty,-c^*)\big)>0$, which is again impossible. Finally, one
has shown that
$$\mu\big((-\infty,-c^*]\cup\{\infty\}\big)=0,$$
that is $\mu$ is supported in $[c^*,+\infty)$.\par
In order to prove that the support of $\mu$ is bounded, it remains to show
that its least upper bound, denoted by $m$, is less than $+\infty$.
The first inequality in (\ref{umu}) implies that, for all
$\gamma\in[c^*,+\infty)$ and $t\in\R$, 
$$u_{\mu}(t,\gamma t)\ge M^{-1}\int_{(\gamma,+\infty)}\varphi_c\big(\gamma
t-ct-c\ln M\big)\,d\mu(c).$$
Whence $\liminf_{t\to+\infty}u_{\mu}(t,\gamma t)\ge
M^{-1}\mu\big((\gamma,+\infty)\big)$ by Lebesgue's dominated convergence
theorem. If $\gamma<m$ then the latter term is positive and thus~(\ref{gtf})
yields $\liminf_{t\to+\infty}\big(X(t)-\gamma t\big)>-\infty$, which in
turn implies
$\liminf_{t\to+\infty}X(t)/t\geq\gamma$. Since this is true for any $\gamma<m$,
but on the other hand we know that $X$ satisfies~\eqref{speedbound} (the 
property~\eqref{speedbound} follows from Proposition~\ref{pro1} which is proved 
independently in Section~\ref{sec44}), we infer that $m<+\infty$. The proof of the 
necessity condition of Theorem~\ref{thsupport} is complete.\hfill$\Box$ 


\subsection{Proof of the sufficiency condition}\label{sec22}

The proof is based on the following proposition, which provides uniform lower
and upper bounds of any solution $u_{\mu}$ of~\eq{P=f} on left and right of
its level sets, when the measure~$\mu\in\mathcal{M}$ is compactly supported in
$[c^*,+\infty)=[2\sqrt{f'(0)},+\infty)$. These bounds say that the~$u_{\mu}$ is
steeper than the standard front associated with any speed larger than its 
support, in the sense of~\eqref{umusteeper3} below,
and they imply especially that the $x$-derivative of $u_{\mu}$ is bounded away
from $0$ along any of its level sets. The following proposition has its own
interest and it will actually be used again in the proof of Theorem~\ref{thspeeds} in Section~\ref{sec3}.

\begin{proposition}\label{pro3}
Under the assumptions~\eqref{hyp-f} and the notations of Section~$\ref{intro}$, let $\mu$ be any measure in~$\mathcal{M}$ that is supported in~$[c^*,\gamma]$ for some $\gamma\in[c^*,+\infty)$, and let $0<u_{\mu}<1$ be the solution of~\eq{P=f} that is associated to the measure $\mu$. Then, for every $(t,X)\in\R^2$,
\be\label{umusteeper3}\left\{\baa{ll}
u_{\mu}(t,X+x)\ge\varphi_{\gamma}\big(\varphi_{\gamma}^{-1}(u_{\mu}(t,X))+x\big) & \hbox{for all }x\le 0,\vspace{3pt}\\
u_{\mu}(t,X+x)\le\varphi_{\gamma}\big(\varphi_{\gamma}^{-1}(u_{\mu}(t,X))+x\big) & \hbox{for all }x\ge 0,\eaa\right.
\ee
where $\varphi_{\gamma}^{-1}:(0,1)\to\R$ denotes the reciprocal of the function $\varphi_{\gamma}$.
\end{proposition}

\begin{remark}{\rm 
Under the assumptions of Proposition~$\ref{pro3}$, it follows from~\eqref{umusteeper3} that the spatial derivative of $u_{\mu}$ at any given point $(t,X)\in\R^2$ is bounded from above by a nega\-tive real number that only depends on the value of $u_{\mu}$ at this point, in the sense that $(u_{\mu})_x(t,X)\le\varphi_{\gamma}'\big(\varphi_{\gamma}^{-1}(u_{\mu}(t,X))\big)<0$ for all $(t,X)\in\R^2$.}
\end{remark}

Postponing the proof of Proposition~\ref{pro3}, let us first complete the proof
of the sufficiency condition of Theorem~\ref{thsupport}. Namely, let
$\mu$ be any measure in $\mathcal{M}$ which is supported in $[c^*,\gamma]$ for
some~$\gamma\in[c^*,+\infty)$ and let us show that the function $u_{\mu}$ is a
transition front connec\-ting~$0$ and $1$ for~\eq{P=f}.\par
As already recalled in the general properties before the statement of Theorem~\ref{thsupport}, one knows that the function $u_{\mu}$ ranges in $(0,1)$ and is nonincreasing with respect to~$x$. Furthermore, the inequalities~(\ref{umu}) read in this case
$$\baa{l}
\max\Big(\varphi_{c^*}\big(x-c^*t-c^*\ln\mu(c^*)\big),\displaystyle M^{-1}\!\!\int_{(c^*,\gamma]}\!\!\varphi_c\big(x\!-\!ct\!-\!c\ln M\big)\,d\mu(c)\Big)\vspace{3pt}\\
\qquad\quad\le\,u_{\mu}(t,x)\,\le\,\varphi_{c^*}\big(x-c^*t-c^*\ln\mu(c^*)\big)+\displaystyle M^{-1}\!\!\int_{(c^*,\gamma]}\!\!e^{-\lambda_c(x-ct-c\ln M)}d\mu(c)\eaa$$
under the usual convention of neglecting the terms involving $\ln0$. Therefore,
$u_{\mu}(t,-\infty)=1$ and $u_{\mu}(t,+\infty)=0$ for every $t\in\R$, from
Lebesgue's dominated convergence theorem. Thus, the function~$v=(u_{\mu})_x$,
which solves $v_t=v_{xx}+f'(u_{\mu})v$ in $\R^2$, cannot be identically equal to
$0$ and is therefore negative in $\R^2$ from the strong maximum principle. In
particular, for every~$m\in(0,1)$ and~$t\in\R$, there is a unique $X_m(t)\in\R$
such that
\be\label{Xmt}
u_{\mu}(t,X_m(t))=m.
\ee\par
We claim that~$u_{\mu}$ is a transition front connecting $0$
and $1$ in the sense of~(\ref{gtf}), with, say,
$X(t)=X_{\varphi_{\gamma}(0)}(t)$. Indeed, applying
\eqref{umusteeper3} in Proposition \ref{pro3} with $X=X(t)$,
we infer that, for every $t\in\R$,
$$\left\{\baa{ll}
u_{\mu}(t,X(t)+x)\ge\varphi_{\gamma}(x)
& \hbox{for all }x\le 0,\vspace{3pt}\\
u_{\mu}(t,X(t)+x)\le\varphi_{\gamma}(x) & \hbox{for all }x\ge 0.\eaa\right.
$$
From this~(\ref{gtf}) follows because $\varphi_{\gamma}(-\infty)=1$ and
$\varphi_{\gamma}(+\infty)=0$. The proof of the sufficiency condition of Theorem~\ref{thsupport} is thereby complete.~\hfill$\Box$


\subsection{Proof of Proposition~\ref{pro3}}\label{sec23}

Proposition~\ref{pro3} states that $u_\mu$ is steeper than any
standard front with a speed lying on the right of the support of $\mu$. It is
natural to expect that $u_\mu$ satisfies this property, since it is constructed
as a superposition of standard fronts with speeds contained in the support of
$\mu$, and it is known, for instance from \cite{aw}, that, if $c<\gamma$,
then $\vp_c$ is steeper than $\vp_\gamma$ in the sense that \eqref{umusteeper3}
holds with $u_\mu$ replaced by $\vp_c$. However, the fact that this property is
preserved after the averaging process in the definition of $u_\mu^n$ is higly
non-trivial. In order to prove it, we first show some bounds similar 
to~(\ref{umusteeper3}) with any speed
$\gamma+\epsilon$ larger than~$\gamma$ and with an approximating sequence
$(u_{\mu_k})_{k\in\N}$. The approximated bounds will actually be uniform with
respect to $\epsilon\in(0,1)$ and $k\in\N$ and we will finally pass to the limit
as $\epsilon\to0^+$ and~$k\to+\infty$. The proof is based on an intersection
number argument and on some uniform decay estimates of the fronts~$\varphi_c$.
  
\subsubsection*{Step 1: uniform decay estimates of the traveling fronts with speeds $c\in[c^*,\gamma]$}

Let us first fix any $\epsilon>0$ ($\epsilon$ will be fixed till almost the end
of the proof of Proposition~\ref{pro3}) and call, in accord with
~(\ref{lambdac}), 
\be\label{deflambda}
\lambda_{\gamma+\epsilon}=\frac{\gamma+\epsilon-\sqrt{(\gamma+\epsilon)^2-4f'(0)
}}{2},
\ee
Notice
that $0<\lambda_{\gamma+\epsilon}<\lambda_{\gamma}\le\lambda_c\le\lambda_{c^*}$
for all~$c\in[c^*,\gamma]$. The introduction of the parameter $\epsilon>0$ is
used for the following auxiliary result. 

\begin{lemma}
There is $A\in\R$ such that
\be\label{claimphic}
\varphi_c'(x)<-\lambda_{\gamma+\epsilon}\,\varphi_c(x)\ \hbox{ for all }x\ge
A\hbox{ and for all }c\in[c^*,\gamma]. 
\ee
\end{lemma}

\noindent{\bf{Proof.}} Assume that there is no real number $A$ such that~(\ref{claimphic}) holds. Then there are two sequences $(x_p)_{p\in\N}$ in $\R$ and $(c_p)_{p\in\N}$ in $[c^*,\gamma]$ such that $x_p\to+\infty$ as~$p\to+\infty$ and
\be\label{varphicp}
\varphi_{c_p}'(x_p)\ge-\lambda_{\gamma+\epsilon}\,\varphi_{c_p}(x_p)\ \hbox{ for all }p\in\N.
\ee
As already recalled in the introduction, one knows that $\varphi_c(x)\le
e^{-\lambda_cx}$ for all $c>c^*$ and for all $x\in\R$ (and
$\varphi_{c^*}(+\infty)=0$). Furthermore,
$0<\lambda_{\gamma}\le\lambda_c\le\lambda_{c^*}$ for all $c\in[c^*,\gamma]$.
Therefore,~$\varphi_{c_p}(x_p)\to0^+$ as $p\to+\infty$. Define, for every
$p\in\N$ and $x\in\R$, 
$$\psi_p(x)=\frac{\varphi_{c_p}(x+x_p)}{\varphi_{c_p}(x_p)}\ (>0).$$
Since each function $\varphi_c$ is positive and satisfies
$\varphi_c''+c\varphi_c'+f(\varphi_c)=0$ in $\R$ and since the sequence
$(c_p)_{p\in\N}$ is bounded, it follows from Harnack inequality that, for every
$R>0$, the sequence $(\psi_p)_{p\in\N}$ is bounded in $L^{\infty}([-R,R])$. On
the other hand, each function $\psi_p$ satisfies
$$\psi_p''+c_p\psi_p'+\frac{f\big(\varphi_{c_p}(x_p)\,\psi_p\big)}{\varphi_{c_p}
(x_p)}=0\ \hbox{ in }\R.$$
Standard elliptic estimates yield the existence of a function $\psi\in C^2(\R)$ and a real number~$c\in[c^*,\gamma]$ such that, up to extraction of a subsequence, $\psi_p\to\psi$ in $C^2_{loc}(\R)$ and $c_p\to c$ as~$p\to+\infty$. The function $\psi$ is nonnegative and obeys
$$\psi''+c\psi'+f'(0)\psi=0\ \hbox{ in }\R.$$
Moreover, $\psi(0)=1$ (whence $\psi$ is positive in $\R$ from the strong maximum principle) and
\be\label{psilambda}
\psi'(0)\ge-\lambda_{\gamma+\epsilon},
\ee
by~(\ref{varphicp}). On the one hand, if $c>c^*=2\sqrt{f'(0)}$, then
$\psi(x)=C_1e^{-\lambda_cx}+C_2e^{-\tilde{\lambda}_cx}$ for all~$x\in\R$, for
some constants $C_1$ and $C_2\in\R$, where
$$\tilde{\lambda}_c=\frac{c+\sqrt{c^2-4f'(0)}}{2}>\frac{c-\sqrt{c^2-4f'(0)}}{2}
=\lambda_c.$$
In this case, $1=\psi(0)=C_1+C_2$ and $C_2$ is nonnegative since $\psi(x)$ is positive in $\R$ and in particular $\psi(x)$ remains positive as $x\to-\infty$. Therefore,
$$\psi'(0)=-\lambda_cC_1-\tilde{\lambda}_cC_2\le-\lambda_cC_1-\lambda_cC_2=-\lambda_c\le-\lambda_{\gamma}<-\lambda_{\gamma+\epsilon}$$
from the choice of $\lambda_{\gamma+\epsilon}$ in~(\ref{deflambda}). This
contradicts~(\ref{psilambda}). On the other hand, if~$c=c^*=2\sqrt{f'(0)}$, then
$\psi(x)=(Cx+1)\,e^{-\lambda_{c^*}x}$ for all $x\in\R$, for some constant
$C\in\R$, with $\lambda_{c^*}=\sqrt{f'(0)}$. The positivity of $\psi$ in $\R$
and in particular as~$x\to-\infty$ implies that $C$ is nonpositive. Therefore,
$\psi'(0)=C-\lambda_{c^*}\le-\lambda_{c^*}<-\lambda_{\gamma+\epsilon}$, and one
has again reached a contradiction. The proof of the claim~(\ref{claimphic}) is
thereby complete.~\hfill$\Box$

\subsubsection*{Step 2: approximation of $u_{\mu}$ by $u_{\mu_k}$}

Now, we approximate $u_{\mu}$ by some functions $u_{\mu_k}$. Namely, when
$\mu(c^*)=0$, we consider any fixed sequence $(\sigma_k)_{k\in\N}$ of positive
real numbers converging to $0^+$, and we define
$$\mu_k=\mu+\sigma_k\delta_{c^*}$$
for every $k\in\N$. 
When $\mu(c^*)>0$, we simply set $\sigma_k=0$ and $\mu_k=\mu$. 
The advantage of adding a positive weight to the point $c^*$ if $\mu(c^*)=0$ is
that it forces the associated function~$u_{\mu_k}$ to have an interface close
to $c^*t$ for large negative $t$. This property, that can be viewed as a
consequence of Theorem \ref{thspeeds} part (ii), will be crucial in the
sequel. Notice that the measures~$\mu_k$ belong to $\mathcal{M}$. Furthermore, 
with the definitions~(\ref{defumun}), for every $k\in\N$ and~$n\in\N$, there holds
\be\label{defumukn}\baa{rcl}
u_{\mu}^n(-n,x) & \!\!\!\!=\!\!\!\! & \max\Big(\varphi_{c^*}\big(x\!+\!c^*n\!-\!c^*\ln\mu(c^*)\big),\displaystyle M^{-1}\!\!\!\int_{(c^*,\gamma]}\!\!\varphi_c\big(x\!+\!cn\!-\!c\ln M\big)d\mu(c)\Big)\vspace{3pt}\\
& \!\!\!\!\le\!\!\!\! & \max\!\Big(\varphi_{c^*}\big(x\!+\!c^*n\!-\!c^*\ln(\mu(c^*)\!+\!\sigma_k)\big),\displaystyle M^{-1}\!\!\!\int_{(c^*,\gamma]}\!\!\!\!\varphi_c\big(x\!+\!cn\!-\!c\ln M\big)d\mu(c)\!\Big)\vspace{3pt}\\
& \!\!\!\!=\!\!\!\! & u^n_{\mu_k}(-n,x)\eaa
\ee
for all $x\in\R$, since $\varphi_{c^*}$ is positive and decreasing, under the
usual convention that the terms involving $\ln 0$ are
absent. For every $k\in\N$, there always holds $u^n_{\mu}(t,x)\le
u^n_{\mu_k}(t,x)$ for all $(t,x)\in[-n,+\infty)\times\R$ from the maximum
principle, and
\be\label{umukmu}
u_{\mu}(t,x)\le u_{\mu_k}(t,x)\ \hbox{ for all }(t,x)\in\R^2
\ee
by passing to the limit $n\to+\infty$.\par
Furthermore, the solutions $u_{\mu_k}$ converge to $u_{\mu}$ locally uniformly in $\R^2$ as $k\to+\infty$.

\begin{lemma}\label{lemumukmu} One has
\be\label{limumukmu}
u_{\mu_k}(t,x)\to u_{\mu}(t,x)\ \hbox{ as }k\to+\infty,\hbox{ locally uniformly in }(t,x)\in\R^2.
\ee
\end{lemma}

\noindent{\bf{Proof.}} If $\mu(c^*)>0$, then $\mu_k=\mu$ and $u_{\mu_k}=u_{\mu}$, so there is nothing to prove. Consider the case~$\mu(c^*)=0$ (whence $M=\mu\big((c^*,+\infty)\big)>0$ in this case). By~(\ref{defumukn}), there holds
$$u_{\mu_k}^n(-n,x)\le\min\big(\varphi_{c^*}(x+c^*n-c^*\ln\sigma_k)+u_{\mu}^n(-n,x),1\big)$$
for every $k$, $n\in\N$ and $x\in\R$. For every $k$ and $n\in\N$, the function $\overline{u}$ defined in $[-n,+\infty)\times\R$ by
$$\overline{u}(t,x)=\min\big(\varphi_{c^*}(x-c^*t-c^*\ln\sigma_k)+u^n_{\mu}(t,x),1\big)$$
is a supersolution for the Cauchy problem satisfied by $u_{\mu_k}^n$ in $[-n,+\infty)\times\R$ since $f(1)=0$ and, when $\overline{u}(t,x)<1$,
$$\baa{rcl}
\overline{u}_t(t,x)-\overline{u}_{xx}(t,x)-f(\overline{u}(t,x)) & = &
f\big(\varphi_{c^*}(x-c^*t-c^*\ln\sigma_k)\big)+f\big(u^n_{\mu}(t,
x)\big)\vspace{3pt}\\
& &
-f\big(\varphi_{c^*}(x-c^*t-c^*\ln\sigma_k)+u^n_{\mu}(t,x)\big)\vspace{3pt},
\eaa$$
which is nonnegative because, for
$a,b\geq0$ with $a+b\le1$, the concavity of $f$ and $f(0)=0$ yield
$$f(a)+f(b)=\int_0^a f'(s)ds+\int_0^b f'(s)ds\geq\int_0^{a+b} f'(s)ds=f(a+b).$$
It follows from the maximum principle that, for
every $k\in\N$ and $n\in\N$,
$$u_{\mu_k}^n(t,x)\le\overline{u}(t,x)=\min\big(\varphi_{c^*}(x\!-\!c^*t\!-\!c^*\ln\sigma_k)\!+\!u_{\mu}^n(t,x),1\big)\hbox{ for all }(t,x)\!\in\![-n,+\infty)\!\times\!\R,$$
whence
$$u_{\mu_k}(t,x)\le\min\big(\varphi_{c^*}(x-c^*t-c^*\ln\sigma_k)+u_{\mu}(t,x),1\big)\ \hbox{ for all }(t,x)\in\R^2$$
by passing to the limit as $n\to+\infty$. Together with~(\ref{umukmu}), one gets that
$$0\le u_{\mu_k}(t,x)-u_{\mu}(t,x)\le\varphi_{c^*}(x-c^*t-c^*\ln\sigma_k)$$
for every $k\in\N$ and $(t,x)\in\R^2$, and the conclusion~(\ref{limumukmu}) follows since $\varphi_{c^*}(+\infty)=0$ and~$\sigma_k\to0^+$. The proof of Lemma~\ref{lemumukmu} is thereby complete.\hfill$\Box$

\subsubsection*{Step 3: the solutions $u_{\mu_k}$ are steeper than $\varphi_{\gamma+\epsilon}$}

The goal of this step is to derive an analogue of Proposition~\ref{pro3} for
the approximated speed~$\gamma+\epsilon$ and the approximated measures $\mu_k$.
Namely:

\begin{lemma}\label{lemsteep}
For every $k\in\N$ and $(t,X)\in\R^2$, there holds
\be\label{steeper3}\left\{\baa{ll}
u_{\mu_k}(t,X+x)\ge\varphi_{\gamma+\epsilon}\big(\varphi_{\gamma+\epsilon}^{-1}(u_{\mu_k}(t,X))+x\big) & \hbox{for all }x\le0,\vspace{3pt}\\
u_{\mu_k}(t,X+x)\le\varphi_{\gamma+\epsilon}\big(\varphi_{\gamma+\epsilon}^{-1}(u_{\mu_k}(t,X))+x\big) & \hbox{for all }x\ge0,\eaa\right.
\ee
whence $(u_{\mu_k})_x(t,X)\le\varphi_{\gamma+\epsilon}'\big(\varphi_{\gamma+\epsilon}^{-1}(u_{\mu_k}(t,X))\big)<0$.
\end{lemma}

\noindent{\bf{Proof.}} Throughout the proof, $k$ denotes any fixed integer. Define
$$B_k=\max\big(A,-c^*\ln(\mu(c^*)+\sigma_k)+\gamma|\ln M|+A\big),$$
where $A\in\R$ is as in~(\ref{claimphic}), and let $u_{k,0}:\R\to(0,1)$ be the continuous function defined by
\be\label{defu0}
u_{k,0}(x)=\left\{\baa{ll}
\varphi_{c^*}(B_k) & \hbox{for }x<0,\vspace{3pt}\\
\varphi_{c^*}(B_k)\,e^{-\lambda_{\gamma+\epsilon}x} & \hbox{for }x\ge0.\eaa\right.
\ee\par
Let us show in this first paragraph that for all $n\in\N$ the functions $u_{\mu_k}^n(-n,\cdot)$ given in~(\ref{defumukn}) are steeper than any translate $u_{k,0}(\cdot+x_0)$ of $u_{k,0}$, in the sense that if the graphs of~$u_{\mu_k}^n(-n,\cdot)$ and~$u_{k,0}(\cdot+x_0)$ intersect, then they have exactly one intersection point and~$u_{\mu_k}^n(-n,\cdot)$ is larger (resp. smaller) than~$u_{k,0}(\cdot+x_0)$ on the left (resp. right) of this intersection point. So, let $n\in\N$ and $x_0\in\R$ be arbitrary and assume that there is $x_1\in\R$ such that
\be\label{defx1}
u_{\mu_k}^n(-n,x_1)=u_{k,0}(x_1+x_0).
\ee
It follows from the definition of the function $u_{\mu_k}^n(-n,\cdot)$ in~(\ref{defumukn}) that $u_{\mu_k}^n(-n,\cdot)$ is continuous and decreasing in $\R$, and that $u_{\mu_k}^n(-n,-\infty)=1$ and $u_{\mu_k}^n(-n,+\infty)=0$. Therefore, there is a unique $\xi_{k,n}\in\R$ such that
\be\label{xin}\left\{\baa{lcll}
u_{\mu_k}^n(-n,\xi_{k,n}) & = & \varphi_{c^*}(B_k), & \vspace{3pt}\\
u_{\mu_k}^n(-n,x) & > & \varphi_{c^*}(B_k)\ (\ge u_{k,0}(x+x_0)) & \hbox{for all }x<\xi_{k,n},\vspace{3pt}\\
u_{\mu_k}^n(-n,x) & < & \varphi_{c^*}(B_k) & \hbox{for all }x>\xi_{k,n}.\eaa\right.
\ee
Since $u_{k,0}\le\varphi_{c^*}(B_k)$ in $\R$ by definition, one infers from~(\ref{defx1}) and~(\ref{xin}) that $x_1\ge\xi_{k,n}$. Furthermore, for all~$x\ge\xi_{k,n}$, one has
$$\varphi_{c^*}\big(x+c^*n-c^*\ln(\mu(c^*)+\sigma_k)\big)\le u^n_{\mu_k}(-n,x)\le\varphi_{c^*}(B_k)$$
by~(\ref{defumukn}) and~(\ref{xin}), whence
$$\baa{rcl}
x & \ge & -c^*n+c^*\ln(\mu(c^*)+\sigma_k)+B_k\vspace{3pt}\\
& \ge & \max\big(\!-c^*n+c^*\ln(\mu(c^*)+\sigma_k)+A,-c^*n+\gamma|\ln M|+A\big)\eaa$$
since $\varphi_{c^*}$ is decreasing. Thus,~(\ref{claimphic}) implies that, for all $x\ge\xi_{k,n}$,
$$\left\{\baa{rcl}
\varphi_{c^*}'\big(x+c^*n-c^*\ln(\mu(c^*)+\sigma_k)\big) & < & -\lambda_{\gamma+\epsilon}\,\varphi_{c^*}\big(x+c^*n-c^*\ln(\mu(c^*)+\sigma_k)\big),\vspace{3pt}\\
\varphi_c'(x+cn-c\ln M) & < & -\lambda_{\gamma+\epsilon}\,\varphi_c(x+cn-c\ln M)\ \hbox{ for all }c^*<c\le\gamma\eaa\right.$$
by~(\ref{claimphic}). As a consequence of the definition~(\ref{defumukn}) of~$u_{\mu_k}^n(-n,\cdot)$, it follows that
$$u_{\mu_k}^n(-n,x')<u_{\mu_k}^n(-n,x)\,e^{-\lambda_{\gamma+\epsilon}(x'-x)}\ \hbox{ for all }\xi_{k,n}\le x<x'.$$
On the other hand, by definition of $u_{k,0}$ in~(\ref{defu0}), there holds
$$u_{k,0}(x'+x_0)\ge u_{k,0}(x+x_0)\,e^{-\lambda_{\gamma+\epsilon}(x'-x)}\ \hbox{ for all }x<x'\in\R.$$
From the definition of $x_1$ in~(\ref{defx1}) and the fact that $x_1\ge\xi_{k,n}$, one concludes that $x_1$ is necessarily unique, and that
$$\left\{\baa{ll}
u_{\mu_k}^n(-n,x)>u_{k,0}(x+x_0) & \hbox{for all }x<\xi_{k,n}\ (\hbox{from }(\ref{xin})),\vspace{3pt}\\
u_{\mu_k}^n(-n,x)>u_{k,0}(x+x_0) & \hbox{for all }\xi_{k,n}\le x<x_1,\vspace{3pt}\\
u_{\mu_k}^n(-n,x)<u_{k,0}(x+x_0) & \hbox{for all }x>x_1.\eaa\right.$$
Finally, for every $n\in\N$ and $x_0\in\R$, the function $u_{\mu_k}^n(-n,\cdot)$ is steeper than $u_{k,0}(\cdot+x_0)$ in the above sense.\par
Let now $u_k$ be the solution of the Cauchy problem~\eq{P=f} with initial condition~$u_{k,0}$ defined in~(\ref{defu0}). It follows from~\cite{a,dm,m2}, as in~\cite{dgm,g}, that, for every $x_0\in\R$, $n\in\N$ and~$t>-n$, the function~$u_{\mu_k}^n(t,\cdot)$ is steeper than $u_k(t+n,\cdot+x_0)$: in other words, if $u_{\mu_k}^n(t,X)=u_k(t+n,X+x_0)$ for some $X\in\R$, then
\be\label{steeper}\left\{\baa{ll}
u_{\mu_k}^n(t,X+x)>u_k(t+n,X+x_0+x) & \hbox{for all }x<0,\vspace{3pt}\\
u_{\mu_k}^n(t,X+x)<u_k(t+n,X+x_0+x) & \hbox{for all }x>0.\eaa\right.
\ee
On the other hand, by~(\ref{defu0}) and since
$0<\lambda_{\gamma+\epsilon}<\lambda_{\gamma}\le\lambda_{c^*}$ and
$\gamma+\epsilon=\lambda_{\gamma+\epsilon}+f'(0)/\lambda_{\gamma+\epsilon}$, it
is also known~\cite{k,u} that there exists a real number $\zeta_k$ such that
\be\label{conv}
u_k(s,x)-\varphi_{\gamma+\epsilon}\big(x-(\gamma+\epsilon)s+\zeta_k\big)\to0\ \hbox{ as }s\to+\infty,\hbox{ uniformly in }x\in\R.
\ee\par
Finally, let $(t,X)$ be any fixed pair in $\R^2$. Since~$u_{\mu_k}^n(t,X)\to u_{\mu_k}(t,X)\in(0,1)$ as~$n\to+\infty$, there holds
$$0<u_{\mu_k}^n(t,X)<\frac{1+u_{\mu_k}(t,X)}{2}\ (<1)\ \hbox{ for all }n\hbox{ large enough}.$$
Furthermore, for each $n>\max(0,-t)$, the function $u_k(t+n,\cdot)$ is decreasing and continuous in~$\R$, while~$\sup_{\R}u_k(t+n,\cdot)\to1$ as $n\to+\infty$ by~(\ref{conv}), and $\inf_{\R}u_k(t+n,\cdot)=u_k(t+n,+\infty)=0$ for every $n>\max(0,-t)$. Therefore, for all $n$ large enough, there is a unique $x_n\in\R$ such that
$$u_k(t+n,X+x_n)=u_{\mu_k}^n(t,X).$$
It follows then from~(\ref{steeper}) that, for all $n$ large enough,
$$\left\{\baa{ll}
u_{\mu_k}^n(t,X+x)>u_k(t+n,X+x_n+x) & \hbox{for all }x<0,\vspace{3pt}\\
u_{\mu_k}^n(t,X+x)<u_k(t+n,X+x_n+x) & \hbox{for all }x>0.\eaa\right.
$$
Thus, by~(\ref{conv}), there exists a sequence $(h_n)_{n\in\N}$
converging to $0^+$ such that
\be\label{steeper2}\left\{\baa{ll}
u_{\mu_k}^n(t,X+x)\ge\varphi_{\gamma+\epsilon}
\big(X+x_n+x-(\gamma+\epsilon)(t+n)+\zeta_k\big)-h_n & \hbox{for all
}x\le0,\vspace{3pt}\\
u_{\mu_k}^n(t,X+x)\le\varphi_{\gamma+\epsilon}
\big(X+x_n+x-(\gamma+\epsilon)(t+n)+\zeta_k\big)+h_n
 & \hbox{for all }x\ge0.\eaa\right.
\ee
From this at $x=0$ and the convergence $u_{\mu_k}^n\to u_{\mu_k}$ as $n\to+\infty$ locally uniformly in $\R^2$, we deduce that
$$\varphi_{\gamma+\epsilon}\big(X+x_n-(\gamma+\epsilon)(t+n)+\zeta_k\big)\to
u_{\mu_k}(t,X)\in(0,1)\ \hbox{ as }n\to+\infty,$$
whence
$$X+x_n-(\gamma+\epsilon)(t+n)+\zeta_k\to\varphi_{\gamma+\epsilon}^{-1}(u_{\mu_k
}(t,X))\ \hbox{ as }n\to+\infty.$$
As a consequence, by passing to the limit as $n\to+\infty$ in~(\ref{steeper2}), the conclusion~(\ref{steeper3}) follows. The proof of
Lemma~\ref{lemsteep} is thereby complete.\hfill$\Box$

\subsubsection*{Step 4: passages to the limit as $k\to+\infty$ and $\epsilon\to0^+$ in the inequalities~\eqref{steeper3}}

Let~$(t,X)$ be arbitrary in $\R^2$. Since $u_{\mu_k}\to u_{\mu}$ locally uniformly in $\R^2$ as~$k\to+\infty$ by Lemma~\ref{lemumukmu}, passing to the limit as~$k\to+\infty$ in~(\ref{steeper3}) yields
\be\label{steepereps}\left\{\baa{ll}
u_{\mu}(t,X+x)\ge\varphi_{\gamma+\epsilon}\big(\varphi_{\gamma+\epsilon}^{-1}(u_{\mu}(t,X))+x\big) & \hbox{for all }x\le0,\vspace{3pt}\\
u_{\mu}(t,X+x)\le\varphi_{\gamma+\epsilon}\big(\varphi_{\gamma+\epsilon}^{-1}(u_{\mu}(t,X))+x\big) & \hbox{for all }x\ge0.\eaa\right.
\ee\par
On the other hand, by standard elliptic estimates and the fact that all solutions $0<v<1$ of $v''+cv'+f(v)=0$ in $\R$ with $v(-\infty)=1$ and $v(+\infty)=0$ are translations of the function~$\varphi_c$, it follows easily that, for every $m\in(0,1)$ and for every sequence $(c_p)_{p\in\N}$ in $[c^*,+\infty)$ such that~$c_p\to c\in[c^*,+\infty)$ as $p\to+\infty$, there holds 
$$\varphi_{c_p}\big(\varphi_{c_p}^{-1}(m)+x\big)\to\varphi_c\big(\varphi_c^{-1}(m)+x\big)\ \hbox{ as }p\to+\infty,\hbox{ uniformly with respect to }x\in\R.$$
Therefore, the desired inequalities~\eqref{umusteeper3} follow
from~(\ref{steepereps}) by passing to the limit as $\epsilon\to0^+$. The proof
of Proposition~\ref{pro3} is thereby complete.\hfill$\Box$ 


\section{Proof of Theorem~\ref{thspeeds}}\label{sec3}

In this section, we carry out the proof of Theorem~\ref{thspeeds}. Namely, when
the measure $\mu\in\mathcal{M}$ is compactly supported in
$[c^*,+\infty)=[2\sqrt{f'(0)},+\infty)$, we prove that $u_{\mu}$ has some
asymptotic past and future speeds and we identify these speeds as the leftmost
and rightmost points of the support of $\mu$. Furthermore, we prove the
convergence of $u_{\mu}$ to some well identified profiles along its level sets
as $t\to\pm\infty$.\par
As in the previous section, an important step in the proof will be the following
proposition, which provides uniform lower and upper bounds of any solution
$u_{\mu}$ of~\eq{P=f} on left and right of its level sets, when the
measure~$\mu\in\mathcal{M}$ is supported in $[\gamma,+\infty)$ for some
$\gamma\ge c^*$. Roughly speaking, the following proposition is
the counterpart of Proposition~\ref{pro3}, in the sense that it states
that~$u_{\mu}$ is less steep than the standard front associated with any speed
on the left of its support. 

\begin{proposition}\label{pro4}
Under the assumptions~\eqref{hyp-f} and the notations of Section~$\ref{intro}$, let $\mu$ be any measure in~$\mathcal{M}$ that is supported in~$[\gamma,+\infty)$ for some $\gamma\ge c^*$ and let $0<u_{\mu}<1$ be the solution of~\eq{P=f} that is associated to the measure $\mu$. Then, for every $(t,X)\in\R^2$,
\be\label{umusteeper4}\left\{\baa{ll}
u_{\mu}(t,X+x)\le\varphi_{\gamma}\big(\varphi_{\gamma}^{-1}(u_{\mu}(t,X))+x\big) & \hbox{for all }x\le 0,\vspace{3pt}\\
u_{\mu}(t,X+x)\ge\varphi_{\gamma}\big(\varphi_{\gamma}^{-1}(u_{\mu}(t,X))+x\big) & \hbox{for all }x\ge 0.\eaa\right.
\ee
\end{proposition}


\subsection{The asymptotic past and future speeds and the limiting profiles}

Postponing the proof of Proposition~\ref{pro4} in Section~\ref{sec32}, let us first complete the proof of Theorem~\ref{thspeeds}. Let $\mu\in\mathcal{M}$ be compactly supported in $[c^*,+\infty)$ and let~$c_-\le c_+\in[c^*,+\infty)$ be the leftmost and rightmost points of the support of $\mu$, that is~$\mu$ is supported in $[c_-,c_+]$, and
$$\mu\big([c_-,c_-+\epsilon]\big)>0\ \hbox{ and }\ \mu\big([\max(c_+-\epsilon,c^*),c_+]\big)>0$$
for every~$\epsilon>0$. Denote
$M=\mu\big(X\backslash\{-c^*,c^*\}\big)=\mu\big((c^*,+\infty)\big)$. As already
proved in Section~\ref{sec22}, we know that $(u_{\mu})_x<0$ in $\R^2$ and that
$0<u_{\mu}<1$ is a transition front connecting~$0$ and $1$ for~\eq{P=f}, namely
there is a family $(X(t))_{t\in\R}$ of real numbers such that~(\ref{gtf}) holds
for $u_{\mu}$. The fact that $c_{\pm}$ are the asymptotic past and future speeds
of $u_{\mu}$ will follow from the fact that the fronts with speeds larger than
$c_-$ (resp. smaller than $c_+$) are hidden as $t\to-\infty$ (resp. as
$t\to+\infty$) by the one with speed $c_-$ (resp. $c_+$). The proof of the
convergence to the asymptotic profiles is subtler and uses in particular the
lower and upper bounds stated in Propositions~\ref{pro3} and~\ref{pro4}.

\subsubsection*{Step 1: $c_-$ is the asymptotic past speed of $u_{\mu}$ and $\limsup_{t\to-\infty}|X(t)-c_-t|<+\infty$ if and only if $\mu(c_-)=\mu\big(\{c_-\}\big)>0$}

We will consider three cases, according to the value of $c_-$ and~$\mu(c_-)$.\par
{\it{Case 1: $\mu(c_-)=0$.}} In this case, $\mu(c^*)=0$, $M>0$ and the inequalities~(\ref{umu}) amount to
$$0<M^{-1}\!\!\int_{(c_-,c_+]}\!\!\varphi_c(x\!-\!ct\!-\!c\ln M)\,d\mu(c)\le u_{\mu}(t,x)\le M^{-1}\!\!\int_{(c_-,c_+]}\!\!e^{-\lambda_c(x-ct-c\ln M)}d\mu(c)$$
for all $(t,x)\in\R^2$. Lebesgue's dominated convergence theorem yields~$u_{\mu}(t,c_-t)\to0$ as~$t\to-\infty$, whence $X(t)-c_-t\to-\infty$ as $t\to-\infty$ by~(\ref{gtf}) and~(\ref{infsup}). In particular, $\liminf_{t\to-\infty}X(t)/t\ge c_-$. Furthermore, for every $\epsilon>0$,
$$\baa{rcl}
\displaystyle\mathop{\liminf}_{t\to-\infty}u_{\mu}(t,(c_-+\epsilon)t) & \ge & \displaystyle\mathop{\liminf}_{t\to-\infty}M^{-1}\!\!\int_{(c_-,c_-+\epsilon)}\!\!\varphi_c\big((c_-+\epsilon)t\!-\!ct\!-\!c\ln M\big)\,d\mu(c)\vspace{3pt}\\
& = & M^{-1}\mu\big((c_-,c_-+\epsilon)\big)\,>\,0,\eaa$$
the positivity following from definition of $c_-$ and the assumption $\mu(c_-)=0$. Therefore, $\liminf_{t\to-\infty}\big(X(t)-(c_-+\epsilon)t\big)>-\infty$, by~(\ref{gtf}), whence $\limsup_{t\to-\infty}X(t)/t\le c_-+\epsilon$, for every $\epsilon>0$. Finally, $X(t)/t\to c_-$ as $t\to-\infty$.\par
{\it{Case 2: $\mu(c_-)>0$ and $c_->c^*$.}} In this case, $M>0$ and the inequalities~(\ref{umu}) read
$$\baa{l}
\displaystyle M^{-1}\varphi_{c_-}(x-c_-t-c_-\ln M)\,\mu(c_-)+M^{-1}\!\!\int_{(c_-,c_+]}\!\!\varphi_c(x\!-\!ct\!-\!c\ln M)\,d\mu(c)\vspace{3pt}\\
\qquad\qquad\displaystyle\le u_{\mu}(t,x)\le M^{-1}e^{-\lambda_{c_-}(x-c_-t-c_-\ln M)}\mu(c_-)+M^{-1}\!\!\int_{(c_-,c_+]}\!\!e^{-\lambda_c(x-ct-c\ln M)}d\mu(c).\eaa$$
Therefore, Lebesgue's dominated convergence theorem yields
$$\baa{rcl}
0\ <\ M^{-1}\varphi_{c_-}(1)\,\mu(c_-) & \le & \displaystyle\mathop{\liminf}_{t\to-\infty}u_{\mu}(t,c_-t+c_-\ln M+1)\vspace{3pt}\\
& \le & \displaystyle\mathop{\limsup}_{t\to-\infty}u_{\mu}(t,c_-t+c_-\ln M+1)\ \le\ M^{-1}e^{-\lambda_{c_-}}\mu(c_-)\ <\ 1,\eaa$$
whence $\limsup_{t\to-\infty}|X(t)-c_-t|<+\infty$ by~(\ref{gtf}) (in particular, $X(t)/t\to c_-$ as $t\to-\infty$).\par
{\it{Case 3: $\mu(c_-)>0$ and $c_-=c^*$.}} In this case, the inequalities~(\ref{umu}) amount to
\be\label{ineqc*c-}\baa{l}
\displaystyle0<\max\Big(\varphi_{c^*}(x-c^*t-c^*\ln\mu(c^*)),M^{-1}\!\!\int_{(c_-,c_+]}\!\!\varphi_c(x\!-\!ct\!-\!c\ln M)\,d\mu(c)\Big)\vspace{3pt}\\
\qquad\quad\displaystyle\le u_{\mu}(t,x)\le\varphi_{c^*}(x-c^*t-c^*\ln\mu(c^*))+M^{-1}\!\!\int_{(c_-,c_+]}\!\!e^{-\lambda_c(x-ct-c\ln M)}d\mu(c),\eaa
\ee
under the convention that the terms involving $M$ are not present if $M=0$ (which is the case if $c_+=c_-=c^*$). It follows from Lebesgue's dominated convergence theorem that~$u_{\mu}(t,c^*t)\to\varphi_{c^*}(-c^*\ln\mu(c^*))\in(0,1)$ as $t\to-\infty$. Therefore,~$\limsup_{t\to-\infty}|X(t)-c^*t|<+\infty$ by~(\ref{gtf}). In particular, $X(t)/t\to c^*=c_-$ as $t\to-\infty$.

\subsubsection*{Step 2: $c_+$ is the asymptotic future speed of $u_{\mu}$ and $\limsup_{t\to-\infty}\big|X(t)-c_+t\big|<+\infty$ if and only if $\mu(c_+)=\mu\big(\{c_+\}\big)>0$}

Namely, one can prove, as in Step~1, that~$X(t)/t\to c_+$ as $t\to+\infty$. Furthermore, $\limsup_{t\to+\infty}|X(t)-c_+t|<+\infty$ if $\mu(c_+)>0$, whereas~$X(t)-c_+t\to-\infty$ as $t\to+\infty$ if $\mu(c_+)=0$. This properties easily follow from the inequalities~(\ref{umu}) and the analysis of the behavior of $u_{\mu}(t,c_+t)$ or $u_{\mu}(t,(c_+-\epsilon)t)$ as $t\to+\infty$.\par
Therefore, at this stage, parts~(i) and~(ii) of Theorem~\ref{thspeeds} have been proved.

\subsubsection*{Step 3: convergence to the profile $\varphi_{c_-}$ as $t\to-\infty$}

The strategy consists in passing to the limit as $t\to-\infty$ along the positions $X(t)$, that is in some functions~$u_{\mu}(t+t_p,X(t_p)+x)$ with $t_p\to-\infty$ as $p\to+\infty$. We shall identify the limit as a solution~$u_{\mu_{\infty}}$ and, by comparisons with approximated solutions, we shall prove that the measure $\mu_{\infty}$ is a multiple of the Dirac mass at the point $c_-$, meaning that $u_{\mu_{\infty}}$ is a standard traveling front with speed $c_-$.\par
We begin with two lemmas on some properties of solutions and limits of solutions of the type $v=u_{\tilde{\mu}}$ when $\tilde{\mu}\in\mathcal{M}$ is compactly supported in $[c^*,+\infty)$ (in this case, we recall that Theorem~\ref{thsupport} implies that $u_{\tilde{\mu}}$ is a transition front connecting $0$ and $1$).

\begin{lemma}\label{lemtauc}
Let $\tilde{\mu}$ be a measure in $\mathcal{M}$ whose support is included in
$[a,b]\subset[c^*,+\infty)$. Let~$(\tilde{X}(t))_{t\in\R}$ be a family such
that~\eqref{gtf} holds for the transition front $u_{\tilde{\mu}}$ connecting~$0$
and~$1$ for~\eq{P=f}. Then, for every $c<a$, there is $\tau_c>0$ such that 
$$\frac{\tilde{X}(t+\tau)-\tilde{X}(t)}{\tau}>c\ \hbox{ for all }t\in\R\hbox{
and }\tau\ge\tau_c.$$ 
\end{lemma}

\noindent{\bf{Proof.}} By~\eqref{gtf} and~\eqref{infsup}, there exists a constant $m\in(0,1)$ such that 
$$u_{\tilde{\mu}}(t,\tilde{X}(t)+x)\geq m\ \hbox{ for all }t\in\R\hbox{ and }x\le0.$$ 
Thus, for any given $t\in\R$, Proposition~\ref{pro4} yields
$$u_{\tilde{\mu}}(t,\tilde{X}(t)\!+\!x)\ge\left\{\baa{rl} m\ge m\,\varphi_a\big(\varphi_a^{-1}(u_{\tilde{\mu}}(t,\tilde{X}(t)))\!+\!x\big) & \!\!\hbox{for all }x\le0,\vspace{3pt}\\
\varphi_a\big(\varphi_a^{-1}(u_{\tilde{\mu}}(t,\tilde{X}(t)))\!+\!x\big)\ge m\,\varphi_a\big(\varphi_a^{-1}(u_{\tilde{\mu}}(t,\tilde{X}(t)))\!+\!x\big) & \!\!\hbox{for all }x\ge0.\eaa\right.$$
The function $m\varphi_a(x-at)$ is a subsolution of~\eq{P=f}, because
$f(0)=0$ together with the concavity of $f$ imply that $f(ms)\geq mf(s)$ 
for any $s\geq0$. As a consequence of the comparison principle, we then
derive
$$u_{\tilde{\mu}}\big(t+\tau,\tilde{X}(t)+x\big)\ge m\,\varphi_a\big(\varphi_a^{-1}
(u_{\tilde{\mu}}(t,\tilde{X}(t)))+x-a\tau\big)\ \hbox{\ for all }
\tau>0,\ x\in\R.$$
In particular, taking $x=a\tau$ yields
$$u_{\tilde{\mu}}\big(t+\tau,\tilde{X}(t)+a\tau\big)\ge
m\,u_{\tilde{\mu}}(t,\tilde{X}(t))\geq m^2\ \hbox{\ for all }
\tau>0.$$
It follows from \eqref{gtf} that there exists $A$, independent of $t$ and
$\tau$, such that 
$$\tilde{X}(t)+a\tau-\tilde{X}(t+\tau)\leq A.$$
The statement of the lemma is an immediate consequence of this inequality.\hfill$\Box$

\begin{lemma}\label{lemmuinfty}
Let $\tilde{\mu}$ be a measure in $\mathcal{M}$ whose support is included in $[a,b]\subset(c^*,+\infty)$. Let~$(t_p,x_p)_{p\in\N}$ be any sequence in $\R^2$ and assume that the functions $v_p(t,x)=u_{\tilde{\mu}}(t_p+t,x_p+x)$ converge in $C^{1,2}_{loc}(\R^2)$ to a solution~$0<v_{\infty}<1$ of~\eq{P=f} in $\R^2$. Then
$$v_{\infty}=u_{\tilde{\mu}_{\infty}},$$
with a measure $\tilde{\mu}_{\infty}\in\mathcal{M}$ which is supported in $[a,b]$.
\end{lemma}

\noindent{\bf{Proof.}} Let first $c$ and $c'$ be any two real numbers such that $c^*<c'<c<a$ and let $\epsilon$ be any positive real number. Since~$u_{\tilde{\mu}}$ is a transition front connecting $0$ and~$1$ for~\eq{P=f}, there is a family~$(\tilde{X}(t))_{t\in\R}$ of real numbers such that~(\ref{gtf}) holds with~$u_{\tilde{\mu}}$. Since~$u_{\tilde{\mu}}(t_p,x_p)=v_p(0,0)\to v_{\infty}(0,0)\in(0,1)$ as $p\to+\infty$, it follows from~(\ref{gtf}) that there is a real number $\tilde{A}\ge0$ such that
$$\big|x_p-\tilde{X}(t_p)\big|\le\tilde{A}\ \hbox{ for all }p\in\N.$$
Furthermore, again from~(\ref{gtf}), there is $A\ge0$ such that
$$u_{\tilde{\mu}}(t,y)\le\epsilon\ \hbox{ for all }t\in\R\hbox{ and }y\ge\tilde{X}(t)+A.$$
On the other hand, Lemma~\ref{lemtauc} yields the existence of $\tau_c>0$ such that
$$\frac{\tilde{X}(t+\tau)-\tilde{X}(t)}{\tau}\ge c\ \hbox{ for all }t\in\R\hbox{ and }\tau\ge\tau_c.$$
Without loss of generality, one can assume that $\tau_c\ge(A+\tilde{A})/(c-c')$. Let now $(t,x)$ be any pair in $\R^2$ such that $t\ge\tau_c$ and $x\ge-c't$. For every $p\in\N$, there holds $\tilde{X}(t_p)-\tilde{X}(t_p-t)\ge ct$, whence
$$\baa{rcl}
x_p+x\ge\tilde{X}(t_p)-\tilde{A}-c't\ge\tilde{X}(t_p-t)+ct-c't-\tilde{A} & \!\!\ge\!\! & \tilde{X}(t_p-t)+(c-c')\tau_c-\tilde{A}\vspace{3pt}\\
& \!\!\ge\!\! & \tilde{X}(t_p-t)+A\eaa$$
and thus $u_{\tilde{\mu}}(t_p-t,x_p+x)\le\epsilon$. This means that
$v_p(-t,x)\le\epsilon$, whence $v_{\infty}(-t,x)\le\epsilon$ by passing to the
limit as $p\to+\infty$. In other words,
$\sup_{x\ge-c't}v_{\infty}(-t,x)\le\epsilon$ for all $t\ge\tau_c$. Thus,
since~$v_{\infty}>0$ in $\R^2$ and $\epsilon>0$ was arbitrary, one gets that
\be\label{vinfty}
\sup_{x\ge-c'|t|}v_{\infty}(t,x)\to0\ \hbox{ as }t\to-\infty.
\ee
Now, since $c'$ can be arbitrary in $(c^*,a)$, Theorem~1.4 of~\cite{hn2} implies that $v_{\infty}$ is of the type
$$v_{\infty}=u_{\tilde{\mu}_{\infty}}$$
for some measure $\tilde{\mu}_{\infty}\in\mathcal{M}$ that is supported in $(-\infty,-a]\cup[a,+\infty)\cup\{\infty\}$. Furthermore, Proposition~\ref{pro3} implies in particular that, for every $x\ge0$ and $p\in\N$,
$$v_p(0,x)=u_{\tilde{\mu}}(t_p,x_p+x)\le\varphi_b\big(\varphi_b^{-1}(u_{\tilde{\mu}}(t_p,x_p))+x\big)=\varphi_b\big(\varphi_b^{-1}(v_p(0,0))+x\big),$$
whence
\be\label{vinfty2}
u_{\tilde{\mu}_{\infty}}(0,x)=v_{\infty}(0,x)\le\varphi_b\big(\varphi_b^{-1}(v_{\infty}(0,0))+x\big)\ \hbox{ for all }x\ge0
\ee
and $v_{\infty}(0,+\infty)=0$. Therefore, as in the proof of the necessity
condition of Theorem~\ref{thsupport} in Section~\ref{sec21}, it follows that the
measure $\tilde{\mu}_{\infty}$ is supported in~$[a,+\infty)$ (otherwise,
$\liminf_{x\to+\infty}v_{\infty}(t,x)>0$ for every $t\in\R$).\par
Let us finally show that $\tilde{\mu}_{\infty}$ is supported in the interval $[a,b]$. Assume by contradiction that there is $b'\in(b,+\infty)$ such that $\tilde{\mu}_{\infty}\big((b',+\infty)\big)>0$. It follows from the inequalities~(\ref{umu}) applied to $v_{\infty}=u_{\tilde{\mu}_{\infty}}$ that, for all $x\in\R$,
$$e^{\lambda_{b'}x}v_{\infty}(0,x)=e^{\lambda_{b'}x}u_{\tilde{\mu}_{\infty}}(0,x)\ge\tilde{M}_{\infty}^{-1}\int_{(b',+\infty)}\min\big(e^{\lambda_{b'}x}\varphi_c(x-c\ln\tilde{M}_{\infty}),1\big)\,d\tilde{\mu}_{\infty}(c),$$
where $\tilde{M}_{\infty}=\tilde{\mu}_{\infty}\big((c^*,+\infty)\big)>0$. But, for every $c>b'\,(>b\ge a>c^*)$, one has $0<\lambda_c<\lambda_{b'}$, whence $e^{\lambda_{b'}x}\varphi_c(x-c\ln\tilde{M}_{\infty})\to+\infty$ as $x\to+\infty$ by~(\ref{asymvarphi}). Therefore, Lebesgue's dominated convergence theorem (in order to apply $1$, we took the minimum with $1$ in the above displayed formula) yields
$$\liminf_{x\to+\infty}\big(e^{\lambda_{b'}x}v_{\infty}(0,x)\big)\ge\tilde{M}_{\infty}^{-1}\tilde{\mu}_{\infty}\big((b',+\infty)\big)>0,$$
whence $e^{\lambda_bx}v_{\infty}(0,x)\to+\infty$ as $x\to+\infty$ since $0<\lambda_{b'}<\lambda_b$. But the inequality~(\ref{vinfty2}) implies actually that $v_{\infty}(0,x)\le Ce^{-\lambda_bx}$ for all $x\ge0$, for some constant $C>0$. One has then reached a contradiction. Finally,~$\tilde{\mu}_{\infty}\big((b',+\infty)\big)=0$ for all $b'>b$ and $\tilde{\mu}_{\infty}$ is supported in~$[a,b]$ (in particular, by Theorem~\ref{thsupport}, the function $v_{\infty}=u_{\tilde{\mu}_{\infty}}$ is a transition front connecting~$0$ and~$1$). The proof of Lemma~\ref{lemmuinfty} is thereby complete.\hfill$\Box$\break

Let us now come back to our transition front $u_{\mu}$ with $\mu\in\mathcal{M}$ being supported in $[c_-,c_+]\subset[c^*,+\infty)$ and $c_{\pm}$ being the leftmost and rightmost points of the support of $\mu$. We assume that either $c_->c^*$, or $c_-=c^*$ and $\mu(c^*)>0$. In order to show that
\be\label{convc-}
u_{\mu}(t,X(t)+\xi(t)+\cdot)\to\varphi_{c_-}\ \hbox{ in }C^2(\R)\ \hbox{ as }t\to-\infty,
\ee
for some bounded function $\xi:(-\infty,0)\to\R$, as well as~(\ref{convc*}) if $c_-=c^*$ and $\mu(c^*)>0$, we are going to pass to the limit as $t\to-\infty$ along a given level set of $u_{\mu}$.\par
First of all, if $c_-=c^*$ and $\mu(c^*)>0$, then the inequalities~(\ref{ineqc*c-}) and Lebesgue's dominated convergence theorem imply that, for every $A\in\R$,
$$\sup_{x\ge A}\big|u_{\mu}(t,x+c^*t)-\varphi_{c^*}(x-c^*\ln\mu(c^*))\big|\to0\ \hbox{ as }t\to-\infty.$$
But $u_{\mu}$ is decreasing in $x$ and less than $1$, while $\varphi_{c^*}(-\infty)=1$. Therefore,
\be\label{conv5}
\sup_{x\in\R}\big|u_{\mu}(t,x+c^*t)-\varphi_{c^*}(x-c^*\ln\mu(c^*))\big|\to0\ \hbox{ as }t\to-\infty.
\ee
Furthermore, the functions $u_{\mu}(t,c^*t+c^*\ln\mu(c^*)+\cdot)$, which converge uniformly to~$\varphi_{c^*}$ as~$t\to-\infty$, actually converge in $C^2(\R)$ by standard parabolic estimates. In other words,~(\ref{convc*}) holds and the proof of part~(iv) of Theorem~\ref{thspeeds} is done.

\begin{remark}{\rm Condition \eqref{convc*} implies \eqref{convc+-} as $t\to-\infty$, 
with $c_-=c^*$. The function~$t\mapsto\xi(t)$ is defined as follows, under the
notations~\eqref{Xmt} in Section~$\ref{sec22}$: set~$m=\varphi_{c^*}(0)\in(0,1)$
and~$\xi(t)=X_m(t)-X(t)$, for all $t\in\R$. Namely,    
$$u_{\mu}(t,X(t)+\xi(t))=u_{\mu}(t,X_m(t))=m=\varphi_{c^*}(0).$$
Notice that the function~$\xi$ is bounded  by~\eqref{gtf}. 
It follows from~\eqref{convc*} and the strict monotonicity of $\varphi_{c^*}$
that $X(t)+\xi(t)-c^*t-c^*\ln\mu(c^*)\to0$ as $t\to-\infty$. Thus,
by~\eqref{convc*} again and by uniform continuity of $\varphi_{c^*}$, one
infers~\eqref{convc+-} as $t\to-\infty$, with $c_-=c^*$.}
\end{remark}

Let us now assume in the sequel, till the end of the proof of Theorem~\ref{thspeeds}, that
$$c_->c^*,$$
and let us complete the proof of part~(iii) of Theorem~\ref{thspeeds}, that is~(\ref{convc+-}).\par
For every $t\in\R$, let $Y(t)=X_{1/2}(t)$ be the unique real number such that
\be\label{Yt}
u_{\mu}(t,Y(t))=\frac12.
\ee
Pick any sequence $(t_p)_{p\in\N}$ of real numbers such that $t_p\to-\infty$ as $p\to-\infty$ and define
\be\label{defvp}
v_p(t,x)=u_{\mu}\big(t_p+t,Y(t_p)+x\big)\ \hbox{ for all }(t,x)\in\R^2.
\ee
From standard parabolic estimates, the functions $v_p$ converge in $C^{1,2}_{loc}(\R^2)$, up to extraction of a subsequence, to a solution $v_{\infty}$ of~\eq{P=f} such that $0\le v_{\infty}\le1$ in $\R^2$ and~$v_{\infty}(0,0)=1/2$, whence $0<v_{\infty}<1$ in $\R^2$ from the strong maximum principle.\par
The key-point of the proof of the convergence~(\ref{convc-}) is the following lemma, from which~(\ref{convc-}) will follow easily by uniqueness of the limit and the fact that the quantities~$X(t)-Y(t)$ are bounded.

\begin{lemma}\label{lemconvc-}
There holds
\be\label{vphic-}
v_{\infty}(t,x)=\varphi_{c_-}\Big(x-c_-t+\varphi_{c_-}^{-1}\Big(\frac12\Big)\Big)\ \hbox{ for all }(t,x)\in\R^2.
\ee
\end{lemma}

\noindent{\bf{Proof.}} As a matter of fact, the proof is easy if $\mu(c_-)>0$, similarly as in the case $c_-=c^*$ and $\mu(c^*)>0$ treated above. The general case $\mu(c_-)=0$ is more difficult due to the fact that the position of $X(t)$ is not just $c_-t+O(1)$ as $t\to-\infty$, since $X(t)-c_-t\to-\infty$. The proof of this lemma is actually done simultaneously in both cases $\mu(c_-)>0$ and $\mu(c_-)=0$.\par
Since $\mu$ is supported in $[c_-,c_+]\subset(c^*,+\infty)$, it follows from Lemma~\ref{lemmuinfty} that $v_{\infty}$ is a transition front connecting $0$ and $1$ of the type
$$v_{\infty}=u_{\mu_{\infty}},$$
for a measure $\mu_{\infty}\in\mathcal{M}$ which is supported in $[c_-,c_+]$. In
order to conclude the proof of Lemma~\ref{lemconvc-}, we shall prove that the
rightmost point of the support of $\mu_{\infty}$ is nothing but $c_-$ itself,
that is $\mu_{\infty}$ is a multiple of the Dirac mass at the point~$c_-$. To do
so, we shall show that $v_{\infty}=u_{\mu_{\infty}}$ can be approximated by some
solutions of~\eq{P=f} which do not move too much faster than $c_-$.\par
If~$c_+=c_-$, then the initial measure $\mu$ is itself a multiple of the Dirac mass at the point~$c_-$; in other words, in this case, there is $\xi\in\R$ such that~$u_{\mu}(t,x)=\varphi_{c_-}(x-c_-t+\xi)$ for all~$(t,x)\in\R^2$, whence~$Y(t)=c_-t-\xi+\varphi_{c_-}^{-1}(1/2)$ for all $t\in\R$ and
$$v_p(t,x)=u_{\mu}\big(t_p+t,Y(t_p)+x\big)=\varphi_{c_-}\Big(x-c_-t+\varphi_{c_-}^{-1}\Big(\frac12\Big)\Big)\ \hbox{ for all }p\in\N\hbox{ and }(t,x)\in\R^2;$$
the desired conclusion~(\ref{vphic-}) is immediate in this case.\par
Therefore, one can focus in the sequel on the case
$$c_-<c_+.$$
Let $c'$ be any real number such that $c_-<c'<c_+$ and denote
\be\label{defmu12}\left\{\baa{ll}
\mu_1=\mu\,\mathbbm{1}_{[c_-,c']}, &
M_1=\mu\big([c_-,c']\big)=\mu_1(X),\vspace{3pt}\\
\mu_2=\mu\,\mathbbm{1}_{(c',c_+]}, &
M_2=\mu\big((c',c_+]\big)=\mu_2(X).\eaa\right.
\ee
Since $c_{\pm}$ are the leftmost and rightmost points of the support of $\mu$,
one has $M_1,M_2>0$ and $M_1+M_2=M=\mu\big([c_-,c_+]\big)$. Remember now
that $u^n_{\mu}(t,x)\to u_{\mu}(t,x)$ as~$n\to+\infty$ for all $(t,x)\in\R^2$,
where $u^n_{\mu}$ solves the Cauchy problem~(\ref{defumun}) in
$[-n,+\infty)\times\R$ with initial condition
$$u^n_{\mu}(-n,x)=M^{-1}\int_{[c_-,c_+]}\varphi_c(x+cn-c\ln M)\,d\mu(c).$$\par
On the one hand, since all functions $\varphi_c$ are decreasing and positive, one has
$$u^n_{\mu}(-n,x)\ge M^{-1}\int_{[c_-,c']}\varphi_c(x+cn-c\ln M_1)\,d\mu(c)=\frac{M_1}{M}\,u^n_{\mu_1}(-n,x)>0$$
for all $n\in\N$ and $x\in\R$. But, for every $n\in\N$, the positive function
$(M_1/M)u^n_{\mu_1}$ is a subsolution of the equation satisfied by $u^n_{\mu}$
in $(-n,+\infty)\times\R$, since $f((M_1/M)s)\ge(M_1/M)f(s)$ for all~$0\le s\le
1$. The maximum principle yields $(M_1/M)u^n_{\mu_1}(t,x)\le u^n_{\mu}(t,x)$ for
all $n\in\N$ and~$(t,x)\in[-n,+\infty)\times\R$, whence 
\be\label{umumu1}
0<\frac{M_1}{M}\,u_{\mu_1}(t,x)\le u_{\mu}(t,x)\ \hbox{ for all }(t,x)\in\R^2,
\ee
by passing to the limit as $n\to+\infty$.\par
On the other hand, calling
\be\label{defxi12}
\xi_1:=c'\ln\frac{M_1}M\ \hbox{  and }\ \xi_2:=c_+\ln\frac{M_2}M,
\ee
one has, for every $n\in\N$ and $x\in\R$,
$$\baa{rcl}
u^n_{\mu}(-n,x) & \le & \displaystyle M^{-1}\int_{[c_-,c']}\varphi_c(x+cn-c\ln M_1+\xi_1)\,d\mu_1(c)\vspace{3pt}\\
& & \displaystyle+M^{-1}\int_{(c',c_+]}\varphi_c(x+cn-c\ln M_2+\xi_2)\,d\mu_2(c)\vspace{3pt}\\
& \le & \displaystyle M_1^{-1}\int_{[c_-,c']}\varphi_c(x+cn-c\ln M_1+\xi_1)\,d\mu_1(c)\vspace{3pt}\\
& & \displaystyle+M_2^{-1}\int_{(c',c_+]}\varphi_c(x+cn-c\ln M_2+\xi_2)\,d\mu_2(c)\vspace{3pt}\\
& = & u^n_{\mu_1}(-n,x+\xi_1)+u^n_{\mu_2}(-n,x+\xi_2).\eaa$$
Since $u^n_{\mu}<1$ in $[-n,+\infty)\times\R$, it follows then as in the proof of Lemma~\ref{lemumukmu} that, for every $n\in\N$, the function $\overline{u}^n$ defined in~$[-n,+\infty)\times\R$ by
$$\overline{u}^n(t,x)=\min\big(u^n_{\mu_1}(t,x+\xi_1)+u^n_{\mu_2}(t,x+\xi_2),1\big)$$
is a supersolution of the equation satisfied by $u^n_{\mu}$. As a consequence,
$$u^n_{\mu}(t,x)\le\min\big(u^n_{\mu_1}(t,x+\xi_1)+u^n_{\mu_2}(t,x+\xi_2),1\big)$$
for every $n\in\N$ and $(t,x)\in[-n,+\infty)\times\R$, whence
\be\label{umumu12}
u_{\mu}(t,x)\le\min\big(u_{\mu_1}(t,x+\xi_1)+u_{\mu_2}(t,x+\xi_2),1\big)\ \hbox{ for all }(t,x)\in\R^2
\ee
by passing to the limit as $n\to+\infty$.\par
Remembering the definition~(\ref{defvp}) of the functions $v_p$, define now, for every $p\in\N$ and~$(t,x)\in\R^2$,
\be\label{defvp12}
v_{p,1}(t,x)=u_{\mu_1}\big(t_p+t,Y(t_p)+x\big)\ \hbox{ and }\ v_{p,2}(t,x)=u_{\mu_2}\big(t_p+t,Y(t_p)+x\big).
\ee
One infers from~(\ref{umumu1}) and~(\ref{umumu12}) that, for every $p\in\N$ and $(t,x)\in\R^2$,
$$\baa{rcl}
\displaystyle\frac{M_1}{M}\,\underbrace{u_{\mu_1}\big(t_p\!+\!t,Y(t_p)\!+\!x\big)}_{=v_{p,1}(t,x)} & \!\!\le\!\! & \underbrace{u_{\mu}\big(t_p\!+\!t,Y(t_p)\!+\!x\big)}_{=v_p(t,x)}\vspace{3pt}\\
& \!\!\le\!\! & \displaystyle\underbrace{u_{\mu_1}\big(t_p\!+\!t,Y(t_p)\!+\!x\!+\!\xi_1\big)}_{=v_{p,1}(t,x+\xi_1)}+\underbrace{u_{\mu_2}\big(t_p\!+\!t,Y(t_p)\!+\!x\!+\!\xi_2\big)}_{=v_{p,2}(t,x+\xi_2)}.\eaa$$
Up to extraction of another subsequence, the functions $v_{p,1}$ and $v_{p,2}$ converge in $C^{1,2}_{loc}(\R^2)$ as~$p\to+\infty$ to two solutions $0\le v_{\infty,1}\le 1$ and $0\le v_{\infty,2}\le 1$ of~\eq{P=f} such that
$$\frac{M_1}{M}\,v_{\infty,1}(t,x)\le v_{\infty}(t,x)\le v_{\infty,1}(t,x+\xi_1)+v_{\infty,2}(t,x+\xi_2)\ \hbox{ for all }(t,x)\in\R^2.$$
But the inequalities~(\ref{umu}) applied to $u_{\mu_2}$ imply that, for every $p\in\N$ and $(t,x)\in\R^2$,
$$0<v_{p,2}(t,x)=u_{\mu_2}(t_p+t,Y(t_p)+x)\le M_2^{-1}\int_{(c',c_+]}e^{-\lambda_c(Y(t_p)+x-c(t_p+t)-c\ln M_2)}d\mu_2(c).$$
Since $c_-<c'$ and $Y(t_p)/t_p\to c_-$ as $p\to+\infty$ (because $X(t)/t\to c_-$ as $t\to-\infty$ from Step~1, and the quantities $Y(t)-X(t)=X_{1/2}(t)-X(t)$ are bounded in $\R$), it follows from Lebesgue's dominated convergence theorem that $v_{p,2}(t,x)\to0$ as $p\to+\infty$ locally uniformly in $(t,x)\in\R^2$, whence $v_{\infty,2}=0$ in $\R^2$ and
\be\label{vinfty3}
\frac{M_1}{M}\,v_{\infty,1}(t,x)\le v_{\infty}(t,x)\le v_{\infty,1}(t,x+\xi_1)\ \hbox{ for all }(t,x)\in\R^2.
\ee
In particular, $v_{\infty,1}>0$ in $\R^2$, while $v_{\infty,1}$ cannot be identically equal to $1$ (otherwise,~$\inf_{\R^2}v_{\infty}\ge M_1/M>0$, which would contradict~(\ref{vinfty}) for instance). Thus, there holds $0<v_{\infty,1}<1$ in $\R^2$ from the strong maximum principle.\par
Furthermore, since $\mu_1$ is supported in~$[c_-,c']\subset(c^*,+\infty)$, it follows from Lemma~\ref{lemmuinfty} that~$v_{\infty,1}$ is a transition front connecting $0$ and $1$ of the type
$$v_{\infty,1}=u_{\mu_{\infty,1}},$$
for a measure $\mu_{\infty,1}\in\mathcal{M}$ which is supported in $[c_-,c']$. From Steps~1 and~2 of the present proof of Theorem~\ref{thspeeds}, the transition front~$v_{\infty}=u_{\mu_{\infty}}$, resp. $v_{\infty,1}=u_{\mu_{\infty,1}}$, has asymptotic past and future speeds which are equal to the leftmost and rightmost points of the support of $\mu_{\infty}$, resp. $\mu_{\infty,1}$. By~(\ref{vinfty3}), the asymptotic past and future speeds of $v_{\infty}$ are equal to those of $v_{\infty,1}$. Therefore, the measures $\mu_{\infty}$ and $\mu_{\infty,1}$ have the same leftmost and rightmost points and, since $\mu_{\infty,1}$ is supported in $[c_-,c']$, so is $\mu_{\infty}$. Finally, since $c'$ was arbitrarily chosen in~$(c_-,c_+)$, one concludes that the measure~$\mu_{\infty}$ is nothing but a multiple of the Dirac mass at~$\{c_-\}$. In other words,~$v_{\infty}(t,x)=u_{\mu_{\infty}}(t,x)=\varphi_{c_-}(x-c_-t+\xi)$ in $\R^2$, for some $\xi\in\R$. Since~$v_{\infty}(0,0)=1/2$, one gets~$\xi=\varphi_{c_-}^{-1}(1/
2)$ and the desired conclusion~(\ref{vphic-}) has been proved. The proof of Lemma~\ref{lemconvc-} is thereby complete.\hfill$\Box$\break

Let us finally complete the proof of the convergence of $u_{\mu}$ along the positions $X(t)$ to the shifted profile $\varphi_{c_-}$ as $t\to-\infty$, that is~(\ref{convc-}), or~(\ref{convc+-}) as $t\to-\infty$, in the case $c_->c^*$. From the uniqueness of the limit of~$u_{\mu}(t_p+t,Y(t_p)+x)$ as $p\to+\infty$ (with $t_p\to-\infty$), one gets that
$$u_{\mu}(t+t',Y(t)+x)\mathop{\longrightarrow}_{t\to-\infty}\varphi_{c_-}\Big(x-c_-t'+\varphi_{c_-}^{-1}\Big(\frac12\Big)\Big),\hbox{ locally uniformly in }(t',x)\in\R^2.$$
Denote $\xi^-(t)=Y(t)-X(t)-\varphi_{c_-}^{-1}(1/2)$ for all $t\in(-\infty,0]$. By~(\ref{gtf}) and $u_{\mu}(t,Y(t))=1/2$, the function $\xi^-$ is then bounded in $(-\infty,0]$. Furthemore,
\be\label{convumu}
u_{\mu}(t,X(t)+\xi^-(t)+x)\to\varphi_{c_-}(x)\ \hbox{ as }t\to-\infty,\hbox{ locally uniformly in }x\in\R.
\ee
Lastly, since $0<u_{\mu}(t,\cdot)<1$ is decreasing in $\R$ for every $t\in\R$, while $\varphi_{c_-}(-\infty)=1$ and~$\varphi_{c_-}(+\infty)=0$, one concludes that the convergence~(\ref{convumu}) is actually uniform in $x\in\R$, and even holds in $C^2(\R)$ from standard parabolic estimates.

\subsubsection*{Step 4: convergence to the profile $\varphi_{c_+}$ along the level sets of $u_{\mu}$ as $t\to+\infty$.}

In this last step, we prove the convergence of $u_{\mu}$ along the positions $X(t)$ as $t\to+\infty$. We recall that $c_->c^*$ here. The proof follows the same lines as for $t\to-\infty$ in Step~3 above. Namely, let $(t_p)_{p\in\N}$ be any sequence of real numbers such that $t_p\to+\infty$ as $p\to+\infty$. Let $(Y(t))_{t\in\R}$ and $(v_p)_{p\in\N}$ be defined as in~(\ref{Yt}) and~(\ref{defvp}). Up to extraction of a subsequence, the functions $v_p$ converge in $C^{1,2}_{loc}(\R^2)$ to a solution $0<v_{\infty}<1$ of~\eq{P=f} such that~$v_{\infty}(0,0)=1/2$. The analogue of Lemma~\ref{lemconvc-} is the following result.

\begin{lemma}\label{lemconvc+}
There holds
\be\label{vphic+}
v_{\infty}(t,x)=\varphi_{c_+}\Big(x-c_+t+\varphi_{c_+}^{-1}\Big(\frac12\Big)\Big)\ \hbox{ for all }(t,x)\in\R^2.
\ee
\end{lemma}

\noindent{\bf{Proof.}} As in the proof of Lemma~\ref{lemconvc-}, we do not distinguish the cases $\mu(c_+)=0$ and~$\mu(c_+)>0$, the latter being actually much easier to handle. From Lemma~\ref{lemmuinfty}, the function $v_{\infty}$ is a transition front connecting~$0$ and~$1$ of the type
$$v_{\infty}=u_{\mu_{\infty}},$$
for a measure $\mu_{\infty}\in\mathcal{M}$ which is supported in $[c_-,c_+]$. This time, our goal is to show that~$\mu_{\infty}$ is a multiple of the Dirac mass at the point $\{c_+\}$. If $c_-=c_+$, then the conclusion~(\ref{vphic+}) is immediate, as in Lemma~\ref{lemconvc-}. Therefore, we focus in the sequel on the case $c_-<c_+$. \par
So, let $c'$ be arbitrary in the interval $(c_-,c_+)$ and let $\mu_1$, $\mu_2$, $M_1$, $M_2$, $\xi_1$ and $\xi_2$ be as in~(\ref{defmu12}) and~(\ref{defxi12}). The same arguments as in the proof of Lemma~\ref{lemconvc-} imply that
$$0<\frac{M_2}{M}\,u_{\mu_2}(t,x)\le u_{\mu}(t,x)\le\min\big(u_{\mu_1}(t,x+\xi_1)+u_{\mu_2}(t,x+\xi_2),1\big)$$
for all $(t,x)\in\R^2$. Let the sequences $(v_{p,1})_{p\in\N}$ and~$(v_{p,2})_{p\in\N}$ be defined as in~(\ref{defvp12}). Up to extraction of a subsequence, they converge locally uniformly in $\R^2$ to some solutions~$0\le v_{\infty,1}\le1$ and $0\le v_{\infty,2}\le1$ of~\eq{P=f} such that
$$\frac{M_2}{M}\,v_{\infty,2}(t,x)\le v_{\infty}(t,x)\le v_{\infty,1}(t,x+\xi_1)+v_{\infty,2}(t,x+\xi_2)$$
for all $(t,x)\in\R^2$. This time, one has
$$0<v_{p,1}(t,x)=u_{\mu_1}(t_p+t,Y(t_p)+x)\le M_1^{-1}\int_{[c_-,c']}e^{-\lambda_c(Y(t_p)+x-c(t_p+t)-c\ln M_1)}d\mu_1(c),$$
but $Y(t_p)/t_p\to c_+$ as $p\to+\infty$ (since $t_p\to+\infty$ as $p\to+\infty$ and $X(t)/t\to c_+$ as $t\to+\infty$). Therefore, Lebesgue's dominated convergence theorem implies that $v_{p,1}(t,x)\to0$ as $p\to+\infty$ locally unifiormly in $(t,x)\in\R^2$, whence $v_{\infty,1}=0$ in $\R^2$ and
\be\label{vinfty4}
\frac{M_2}{M}\,v_{\infty,2}(t,x)\le v_{\infty}(t,x)\le v_{\infty,2}(t,x+\xi_2)\ \hbox{ for all }(t,x)\in\R^2.
\ee
As a consequence, $0<v_{\infty,2}<1$ in $\R^2$. This time, since $\mu_2$ is supported in~$[c',c_+]\subset(c^*,+\infty)$, it follows from Lemma~\ref{lemmuinfty} that $v_{\infty,2}$ is a transition front connecting $0$ and $1$ of the type
$$v_{\infty,2}=u_{\mu_{\infty,2}},$$
for a measure $\mu_{\infty,2}\in\mathcal{M}$ which is supported in $[c',c_+]$. From Steps~1 and~2, the transition front~$v_{\infty}=u_{\mu_{\infty}}$, resp. $v_{\infty,2}=u_{\mu_{\infty,2}}$, has asymptotic past and future speeds which are equal to the leftmost and rightmost points of the support of $\mu_{\infty}$, resp. $\mu_{\infty,2}$. By~(\ref{vinfty4}), the asymptotic past and future speeds of $v_{\infty}$ are equal to those of $v_{\infty,2}$. Therefore, the measures~$\mu_{\infty}$ and~$\mu_{\infty,2}$ have the same leftmost and rightmost points and, since $\mu_{\infty,2}$ is supported in $[c',c_+]$, so is $\mu_{\infty}$. Finally, since $c'$ was arbitrarily chosen in $(c_-,c_+)$, one concludes that the measure~$\mu_{\infty}$ is nothing but a multiple of the Dirac mass at~$\{c_+\}$. The conclusion of Lemma~\ref{lemconvc+} follows immediately.\hfill$\Box$\break

By denoting $\xi^+(t)=Y(t)-X(t)-\varphi_{c_+}^{-1}(1/2)$ for all $t\in(0,+\infty)$, one concludes as at the end of Step~3 that the function $\xi^+$ is bounded in $(0,+\infty)$ and that
$$u_{\mu}(t,X(t)+\xi^+(t)+\cdot)\to\varphi_{c_+}\ \hbox{ in }C^2(\R)\ \hbox{ as }t\to+\infty.$$
Finally,~(\ref{convc+-}) holds with $\xi(t)=\xi^-(t)$ for $t\le0$ and $\xi(t)=\xi^+(t)$ for $t>0$. The proof of Theorem~\ref{thspeeds} is thereby complete.\hfill$\Box$


\subsection{Proof of Proposition~\ref{pro4}}\label{sec32}

Let $\mu\in\mathcal{M}$ be supported in $[\gamma,+\infty)$ for some $\gamma\ge
c^*=2\sqrt{f'(0)}$ and let $0<u_{\mu}<1$ be the solution of~\eq{P=f} associated
to the measure $\mu$. We will prove that $u_{\mu}$ is less steep than the
traveling front $\varphi_{\gamma}$ associated to the speed $\gamma$. As for
Proposition~\ref{pro3}, the proof uses an intersection number argument. But the
proof is a bit simpler here since one does not need to approximate $u_{\mu}$ by
some $u_{\mu_k}$ nor to approximate the speed $\gamma$ by
$\gamma+\epsilon$.\par 
To do so, let us first show that every front $\varphi_c$, which is already known to be such that~$\varphi_c'<0$ in $\R$, does not decay logarithmically faster than its decay rate at~$+\infty$. Namely, let~$c$ be any real number in~$[c^*,+\infty)$. It is known, see e.g.~\cite{aw}, that~$\phi_c(\xi):=\varphi_c'(\xi)/\varphi_c(\xi)\to-\lambda_c$ as~$\xi\to+\infty$. We claim that
\be\label{claimlambda}
-\lambda_c<\frac{\varphi_c'(\xi)}{\varphi_c(\xi)}=\phi_c(\xi)\ (<0)\ \hbox{ for all }\xi\in\R.
\ee
Indeed, if this were not true, then there would exist a real number $\xi_m$ such that $\phi_c(\xi_m)\le\phi_c(\xi)<0$ for all $\xi\in\R$. But the function $\phi_c$ obeys
$$\phi_c''+(2\phi_c+c)\phi_c'+\Big(f'(\varphi_c)-\frac{f(\varphi_c)}{\varphi_c}\Big)\,\phi_c=0\ \hbox{ in }\R,$$
while $f'(\varphi_c)-f(\varphi_c)/\varphi_c\le0$ in $\R$ since $f(s)/s$ is nonincreasing on $(0,1]$. The strong maximum principle would imply that $\phi_c$ is constant, which is clearly impossible since $\phi_c(-\infty)=0$ and~$\phi_c<0$ in $\R$. Therefore, the claim~(\ref{claimlambda}) holds.\par
Let now $v$ be the solution of the Cauchy problem
$$\left\{\baa{rcl}
v_t & = & v_{xx}+f(v),\ \ t>0,\ x\in\R,\vspace{3pt}\\
v_0(x) & = & \left\{\baa{ll}
1 & \hbox{for all }x\le 0,\vspace{3pt}\\
e^{-\lambda_{\gamma}x} & \hbox{for all }x>0.\eaa\right.\eaa\right.$$
It follows from~\cite{u} that
\be\label{conv3}
\sup_{x\in\R}\big|v(t,x)-\varphi_{\gamma}(x-\gamma t+m(t))\big|\to0\ \hbox{ as }t\to+\infty,
\ee
where $m(t)$ can be taken to be $0$ if $\gamma>c^*$ (remember that $\varphi_c(\xi)\sim e^{-\lambda_c\xi}$ as $\xi\to+\infty$ if~$c>c^*$) and $m(t)=o(t)$ as $t\to+\infty$ if $\gamma=c^*$.\par
On the other hand, by definition, the solution $u_{\mu}$ of~\eq{P=f} is given as the limit as $n\to+\infty$ of the solutions~$u_{\mu}^n$ of the Cauchy problems~(\ref{defumun}) in $[-n,+\infty)\times\R$ associated with the Cauchy data
$$u_{\mu}^n(-n,x)=\max\Big(\varphi_{c^*}(x+c^*n-c^*\ln\mu(c^*)),M^{-1}\int_{(c^*,+\infty)}\varphi_c(x+cn-c\ln M)\,d\mu(c)\Big)$$
for all $x\in\R$, under the convention that the terms involving $\ln0$ are
absent. For every $c\ge\gamma\,(\ge c^*)$, there
holds $0<\lambda_c\le\lambda_{\gamma}$. Since the support of $\mu$ is included
in $[\gamma,+\infty)$, the property~(\ref{claimlambda}) implies that, for all
$n\in\N$,
$$u_{\mu}^n(-n,x')>u_{\mu}^n(-n,x)\,e^{-\lambda_{\gamma}(x'-x)}\ \hbox{ for all }x<x'\in\R.$$
Since the continuous function $u_{\mu}^n(-n,\cdot)$ ranges in $(0,1)$, one infers that $v_0(x_0+\cdot)$ is steeper than $u_{\mu}^n(-n,\cdot)$ for every $n\in\N$ and $x_0\in\R$, in the sense that there is a unique $x_1\in\R$ such that
$$\left\{\baa{rcll}
v_0(x_0+x_1) & = & u_{\mu}^n(-n,x_1), & \vspace{3pt}\\
v_0(x_0+x_1+x) & > & u_{\mu}^n(-n,x_1+x) & \hbox{for all }x<0,\vspace{3pt}\\
v_0(x_0+x_1+x) & < & u_{\mu}^n(-n,x_1+x) & \hbox{for all }x>0.\eaa\right.$$
As a consequence, as in the proof of Lemma~\ref{lemsteep} for instance, it follows from an intersection number argument that, for every $n\in\N$, $t>-n$ and $x_0\in\R$, $v(t+n,x_0+\cdot)$ is steeper than~$u_{\mu}^n(t,\cdot)$.\par
Let finally $(t,X)$ be any fixed pair in $\R^2$. For every $n\in\N$ such that $t>-n$, since~$v(t+n,\cdot+X)$ is continuously decreasing and converges to $0$ and $1$ at $\pm\infty$, there is a unique real number~$x_n$ such that $v(t+n,x_n+X)=u_{\mu}^n(t,X)$. Since $v(t+n,x_n+\cdot)$ is steeper than~$u_{\mu}^n(t,\cdot)$, one infers that
$$\left\{\baa{ll}
v(t+n,x_n+X+x)>u_{\mu}^n(t,X+x) & \hbox{for all }x<0,\vspace{3pt}\\
v(t+n,x_n+X+x)<u_{\mu}^n(t,X+x) & \hbox{for all }x>0.\eaa\right.$$
Thus, by~(\ref{conv3}), there exists a sequence $(h_n)_{n\in\N}$
converging to $0^+$ such that
\be\label{steeper5}\left\{\baa{ll}
u_{\mu}^n(t,X+x)\le\varphi_{\gamma}
\big(x_n+X+x-\gamma(t+n)+m(t+n)\big)+h_n & \hbox{for all
}x\le0,\vspace{3pt}\\
u_{\mu}^n(t,X+x)\ge\varphi_{\gamma}
\big(x_n+X+x-\gamma(t+n)+m(t+n)\big)-h_n & \hbox{for all
}x\ge0.\eaa\right.
\ee
From this at $x=0$ and the convergence $u_{\mu}^n\to u_{\mu}$ as $n\to+\infty$ locally uniformly in $\R^2$, we deduce that
$$\varphi_{\gamma}\big(x_n+X-\gamma(t+n)+m(t+n)\big)\to
u_{\mu}(t,X)\in(0,1)\ \hbox{ as }n\to+\infty,$$
whence
$$x_n+X-\gamma(t+n)+m(t+n)\to\varphi_{\gamma}^{-1}(u_{\mu
}(t,X))\ \hbox{ as }n\to+\infty.$$
Passing to the limit as $n\to+\infty$ in~(\ref{steeper5}) leads to the desired conclusion~(\ref{umusteeper4}) and the proof of Proposition~\ref{pro4} is thereby complete.\hfill$\Box$


\subsection{An additional comparison result for any time-global solution $u$ of~\eq{P=f}}

To complete this section, we include an additional comparison result which
has its own inte\-rest. This
result is of the same spirit as Proposition~\ref{pro4} but it holds for any
solution~$0<u<1$ of~\eq{P=f}, even not a transition front or even not a
solution which is a priori of the type $u_{\mu}$. It says that the traveling front $\varphi_{c^*}(x-c^*t)$ is the steepest time-global solution (that is a critical front, in the sense of~\cite{n2}). This result, stated in Proposition~\ref{pro5} below, actually holds
under more general assumptions on~$f$, for which solutions of the type $u_{\mu}$
are not known. Namely, we will consider~$C^1([0,1])$ functions~$f$ satisfying
\be\label{hypg}
f(0)=f(1)=0\ \hbox{ and }\ 0<f(u)\le f'(0)u\hbox{ for all }u\in(0,1).
\ee
For such functions $f$, standard traveling fronts $\varphi_c(x-ct)$
solving~\eq{P=f} and such that
$1=\varphi_c(-\infty)>\varphi_c>\varphi_c(+\infty)=0$ still exist if and only if $c\ge
c^*=2\sqrt{f'(0)}$.

\begin{proposition}\label{pro5}
Let $f:[0,1]\to\R$ be any $C^1$ function satisfying~\eqref{hypg}. Then, for
every solution $0<u<1$ of~\eq{P=f} and for every $(t,X)\in\R^2$, there holds
\be\label{umusteeper6}\left\{\baa{ll}
u(t,X+x)\le\varphi_{c^*}\big(\varphi_{c^*}^{-1}(u(t,X))+x\big) & \hbox{for all
}x\le 0,\vspace{3pt}\\
u(t,X+x)\ge\varphi_{c^*}\big(\varphi_{c^*}^{-1}(u(t,X))+x\big) & \hbox{for all
}x\ge 0.\eaa\right.
\ee
\end{proposition}

\noindent{\bf{Proof.}} In the proof, $0<u<1$ denotes any solution of~\eq{P=f}
and~$(t,X)$ denotes any fixed pair in $\R^2$. The proof will again use an
intersection number argument, by comparing this time~$u$ with the solution of
some Cauchy problems with the steepest possible initial condition. Namely,
let~$v$ be the solution of the Cauchy problem
$$\left\{\baa{rcl}
v_{\tau} & = & v_{xx}+f(v),\ \ \tau>0,\ x\in\R,\vspace{3pt}\\
v(0,x) & = & \left\{\baa{ll} 1 & \hbox{for all }x<0,\vspace{3pt}\\
0 & \hbox{for all }x\ge 0.\eaa\right.\eaa\right.$$
From the maximum principle and the parabolic regularity, the function
$v(\tau,\cdot)$ is decreasing and continuous for every $\tau>0$, while
$v(\tau,-\infty)=1$ and $v(\tau,+\infty)=0$. Therefore, for every~$\tau>0$,
there is a unique $m(\tau)\in\R$ such that
$$v(\tau,m(\tau))=u(t,X)\ (\in(0,1)).$$
It is also known, see e.g.~\cite{b,dgm,hnrr,kpp,l,u},
that~$m(\tau)=c^*\tau-(3/c^*)\ln\tau+O(1)$ as~$\tau\to+\infty$ and that
\be\label{convv}
v(\tau,m(\tau)+x)\to\varphi_{c^*}\big(\varphi_{c^*}^{-1}(u(t,X))+x\big)\ \hbox{
as }t\to+\infty,\hbox{ uniformly in }x\in\R,
\ee
where $\varphi_{c^*}(x-c^*t)$ denotes any traveling front with critical speed
$c^*=2\sqrt{f'(0)}$ for the equation~\eq{P=f}.\par
Now, owing to the definition of $v_0$, for every $s\in(-\infty,t)$, the shifted
step-function~$v(0,m(t-s)+\cdot)$ is steeper than $u(s,X+\cdot)$, in the sense
that
$$\left\{\baa{ll}
u(s,X+x)<1=v(0,m(t-s)+x) & \hbox{for all }x<-m(t-s),\vspace{3pt}\\
u(s,X+x)>0=v(0,m(t-s)+x) & \hbox{for all }x>-m(t-s).\eaa\right.$$
Therefore, it follows from an intersection number argument
that~$v(t-s,m(t-s)+\cdot))$ is steeper than $u(t,X+\cdot)$ and, since
$$v(t-s,m(t-s))=u(t,X),$$
one infers that
$$\left\{\baa{ll}
u(t,X+x)<v(t-s,m(t-s)+x) & \hbox{for all }x<0,\vspace{3pt}\\
u(t,X+x)>v(t-s,m(t-s)+x) & \hbox{for all }x>0.\eaa\right.$$
By passing to the limit as $s\to-\infty$ and using~(\ref{convv}), the
conclusion~\eqref{umusteeper6} follows and the proof of Proposition~\ref{pro5}
is thereby complete.\hfill$\Box$


\section{Proofs of the other main results}\label{sec4}

In this section, we do the proofs of Theorems~\ref{thmeanspeed},~\ref{thliminf},~\ref{thm:decay} and Proposition~\ref{pro1}. Theorems~\ref{thmeanspeed},~\ref{thliminf},~\ref{thm:decay} will chiefly follow from Theorems~\ref{thsupport} and~\ref{thspeeds}, while the proof of Proposition~\ref{pro1} will be done in a more general heterogeneous framework, in Section~\ref{sec44}.


\subsection{Proof of Theorem~\ref{thm:decay}}

The assumption~(\ref{hyphn}) with $c>c^*$ implies, by Theorem~1.4 of~\cite{hn2}, that $u$ is of the type~$u=u_{\mu}$, with a measure $\mu\in\mathcal{M}$ that is supported in $(-\infty,-c]\cup[c,+\infty)\cup\{\infty\}$. In particular, $M=\mu\big(X\backslash\{-c^*,c^*\}\big)=\mu(X)>0$. Two cases may occur:\par
{\it Case 1: $\mu\big((-\infty,-c]\cup\{\infty\}\big)>0$.} In this case, it follows from~(\ref{umu}) and Lebesgue's dominated convergence theorem that
$$\liminf_{x\to+\infty}\underbrace{u(0,x)}_{=u_{\mu}(0,x)}\ge\frac{\mu\big((-\infty,-c]\big)+\theta(\ln M)\,\mu(\infty)}{M}>0.$$
Since $u(0,x)<1$ for all $x\in\R$, it follows then that
$$\frac{\ln u(0,x)}{x}\to0\ \hbox{ as }x\to+\infty.$$\par
{\it Case 2: $\mu\big((-\infty,-c]\cup\{\infty\}\big)=0$.} In this case, $\mu$ is supported in $[c,+\infty)\subset(c^*,+\infty)$ and the inequalities~(\ref{umu}) imply that
$$M^{-1}\int_{[c,+\infty)}\varphi_{c'}(x-c'\ln M)\,d\mu(c')\le\underbrace{u(0,x)}_{=u_{\mu}(0,x)}\le M^{-1}\int_{[c,+\infty)}e^{-\lambda_{c'}(x-c'\ln M)}\,d\mu(c')$$
for all $x\in\R$.\par
We now claim that, for every $\gamma\in(c^*,+\infty)$ such that $\mu\big([\gamma,+\infty)\big)>0$, one has
\be\label{claimlambda2}
\limsup_{x\to+\infty}-\frac{\ln u(0,x)}{x}\le\lambda_{\gamma}=\frac{\gamma-\sqrt{\gamma^2-4f'(0)}}{2}.
\ee
Indeed, for such a $\gamma$ and for any $\lambda'>\lambda_{\gamma}$, we have
$$e^{\lambda'x}u(0,x)\ge
M^{-1}\int_{[\gamma,+\infty)}\min\big(e^{\lambda'x}\varphi_{c'}(x-c'\ln
M),1\big)\,d\mu(c')\ \hbox{ for all }x\in\R.$$
For every $c'\ge\gamma\,(>c^*)$, there holds $\varphi_{c'}(x-c'\ln M)\sim
e^{-\lambda_{c'}(x-c'\ln M)}$ as $x\to+\infty$ with
$0<\lambda_{c'}\le\lambda_{\gamma}<\lambda'$, whence
$$\min\big(e^{\lambda'x}\varphi_{c'}(x-c'\ln M),1\big)\to1\ \hbox{ as }x\to+\infty.$$
Finally, Lebesgue's dominated convergence theorem shows that
$$\liminf_{x\to+\infty}\big(e^{\lambda'x}u(0,x)\big)\ge M^{-1}\,\mu\big([\gamma,+\infty)\big)>0,$$
whence $\limsup_{x\to+\infty}-(\ln u(0,x))/x\le\lambda'$. Since $\lambda'$ can be arbitrary in $(\lambda_{\gamma},+\infty)$, the claim~(\ref{claimlambda2}) follows.\par
{\it Subcase 2.1: the measure $\mu$ is not compactly supported.} Therefore, it follows from~(\ref{claimlambda2}) that, for every $\gamma\in(c^*,+\infty)$, $\limsup_{x\to+\infty}-(\ln u(0,x))/x\le\lambda_{\gamma}=(\gamma-\sqrt{\gamma^2-4f'(0)})/2$, whence $\limsup_{x\to+\infty}-(\ln u(0,x))/x\le0$ by passing to the limit as $\gamma\to+\infty$. But since~$u(0,x)\in(0,1)$ for all $x\in\R$, one concludes in this case that
$$\frac{\ln u(0,x)}{x}\to0\ \hbox{ as }x\to+\infty.$$\par
{\it Subcase 2.2: the measure $\mu$ is compactly supported.} In this case, call $c_+$ the rightmost point of the measure $\mu$ and notice that $c^*<c\le c_+<+\infty$. It follows from Proposition~\ref{pro3} that
$$\underbrace{u(0,x)}_{=u_{\mu}(0,x)}\le\varphi_{c_+}\big(\varphi_{c_+}^{-1}(u(0,0))+x\big)\ \hbox{for all }x\ge 0.$$
Since $\varphi_{c_+}(\xi)\sim e^{-\lambda_{c_+}\xi}$ as $\xi\to+\infty$, one infers that $\limsup_{x\to+\infty}(\ln u(0,x))/x\le-\lambda_{c_+}$, that is $\liminf_{x\to+\infty}(-\ln u(0,x))/x\ge\lambda_{c_+}$. But $\mu\big([\gamma,+\infty)\big)>0$ for every $\gamma\in(c^*,c_+)$, by definition of $c_+$. Thus,~(\ref{claimlambda2}) implies that $\limsup_{x\to+\infty}(-\ln u(0,x))/x\le\lambda_{\gamma}$ for every $\gamma\in(c^*,c_+)$, whence $\limsup_{x\to+\infty}(-\ln u(0,x))/x\le\lambda_{c_+}$. Finally,
\be\label{lambdac+}
-\frac{\ln u(0,x)}{x}\to\lambda_{c_+}\ \hbox{ as }x\to+\infty,
\ee
with $0<\lambda_{c_+}\le\lambda_c<\lambda_{c^*}=\sqrt{f'(0)}$.\par
{\it Conclusion.} To sum up, in all above cases, there holds $u=u_{\mu}$,
$$-\frac{\ln u(0,x)}{x}\to\lambda\in\big[0,\sqrt{f'(0)}\big)\ \hbox{ as
}x\to+\infty$$
and $\lambda>0$ if and only if the measure $\mu$ is compactly supported in $[c^*,+\infty)$. In this case (that is subcase 2.2 with the previous notations), call $c_-$ the leftmost point of the support of $\mu$ (one has $c^*<c\le c_-\le c_+<+\infty$). It follows from Theorems~\ref{thsupport} and~\ref{thspeeds} that $u$ is a transition front connecting $0$ and $1$, that it has asymptotic past and future speeds equal to~$c_{\pm}$, and that~(\ref{convc+-}) holds for some bounded function $\xi$. Furthermore,
$$c_+=\lambda_{c_+}+\frac{f'(0)}{\lambda_{c_+}}=\lambda+\frac{f'(0)}{\lambda}$$
by~(\ref{lambdac+}). On the other hand, since $X(t)/t\to c_-$ as $t\to-\infty$, where $(X(t))_{t\in\R}$ is a family of real numbers such that~(\ref{gtf}) holds, one has $\max_{[-c'|t|,c'|t|]}u(t,\cdot)\to0$ as $t\to-\infty$ for every~$c'\in[0,c_-)$, while~(\ref{infsup}) implies that $\liminf_{t\to-\infty}\max_{[-c'|t|,c'|t|]}u(t,\cdot)>0$ for every $c'>c_-$. As a conclusion,
$$c_-=\sup\Big\{c\ge0,\ \lim_{t\to-\infty}\max_{[-c|t|,c|t|]}u(t,\cdot)=0\Big\}.$$
The proof of Theorem~\ref{thm:decay} is thereby complete.\hfill$\Box$


\subsection{Proof of Theorem~\ref{thliminf}}

Let $u$ be a transition front connecting $0$ and $1$ for~\eq{P=f}, and let $(X(t))_{t\in\R}$ be a family of real numbers such that~(\ref{gtf}) holds. It follows from~(\ref{gtf}) and~(\ref{awspreading}) that $\liminf_{t\to-\infty}X(t)/t\ge c^*$, that is the first inequality in~(\ref{inequalities}).\par
On the other hand, the last inequality in~(\ref{inequalities}) holds by the
consequence \eqref{speedbound} of Proposition~\ref{pro1}.\par
Furthermore, the spreading results of Aronson and Weinberger~\cite{aw} also imply that $\min_{[-ct,ct]}u(t,\cdot)\to1$ as $t\to+\infty$ for every $c\in[0,c^*)$, whence $\liminf_{t\to+\infty}X(t)/t\ge c^*$ by~(\ref{gtf}). In particular, if $c^*=\liminf_{t\to-\infty}X(t)/t$, then formula~(\ref{inequalities}) holds.\par
Lastly, if $c^*<\liminf_{t\to-\infty}X(t)/t$, then~(\ref{gtf}) implies that $\max_{[-c|t|,c|t|]}u(t,\cdot)\to0$ as~$t\to-\infty$ for some $c>c^*$. It follows from Theorem~1.4 of~\cite{hn2} that $u=u_{\mu}$ for some measure~$\mu\in\mathcal{M}$ that is supported in $(-\infty,-c]\cup[c,+\infty)\cup\{\infty\}$. From Theorem~\ref{thsupport}, one infers that~$\mu$ is actually compactly supported in $[c,+\infty)\subset(c^*,+\infty)$ and the conclusion of Theorem~\ref{thliminf} follows from Theorem~\ref{thspeeds}.\hfill$\Box$


\subsection{Proof of Theorem~\ref{thmeanspeed}}

First of all, all speeds $\gamma$ in $[c^*,+\infty)$ are admissible global mean
speeds since every standard traveling front~$\varphi_{\gamma}(x-\gamma t)$ has a
global mean speed equal to $\gamma$. Furthermore, if a
transition front connecting~$0$ and $1$ has a global mean speed equal to
$\gamma$ then necessarily $\gamma\geq c^*$ by the spreading result~\cite{aw}
(in the form of~(\ref{awspreading}), or even by Theorem~\ref{thasymptotic}).
Now, let $u$ be a transition front connecting $0$ and $1$ and having a global
mean speed $\gamma$ such that $\gamma>c^*$, whence $X(t)/t\to\gamma>c^*$ as
$t\to-\infty$ and $\max_{[-c'|t|,+\infty)}u(t,\cdot)\to0$ as $t\to-\infty$, for
every~$c'\in(c^*,\gamma)$. Theorem~1.4 of~\cite{hn2} then implies that $u$ is of
the type 
$$u=u_{\mu},$$
for a measure $\mu\in\mathcal{M}$ whose support is included in $(-\infty,-c']\cup[c',+\infty)\cup\{\infty\}$. On the other hand, since $u=u_{\mu}$ is a transition front connecting~$0$ and~$1$ for~\eq{P=f}, it follows from Theorem~\ref{thsupport} that $\mu$ is compactly supported in $[c^*,+\infty)$. Since $u_{\mu}$ has asymptotic past and future speeds which are both equal to $c$, Theorem~\ref{thspeeds} implies that the leftmost and rightmost points of the support of $\mu$ are both equal to $c$. In other words,~$\mu=m\,\delta_c$ for some~$m>0$ and the conclusion follows immediately.\hfill$\Box$


\subsection{Proof of Proposition~\ref{pro1}}\label{sec44}

Proposition~\ref{pro1} is actually a particular case of Proposition~\ref{pro1bis} below, which holds in a more general framework. More precisely, we consider in this section some heterogeneous reaction-diffusion equations of the type
\be\label{eqbis}
u_t(t,x)=a(t,x)\,u_{xx}(t,x)+b(t,x)\,u_x(t,x)+f(t,x,u(t,x)),\ \ t\in\R,\ x\in\R,
\ee
where $a:\R^2\to\R$ is a bounded continuous function such that $\inf_{\R^2}a>0$, $b:\R^2\to\R$ is a bounded continuous function, and $f:\R^2\times[0,1]\to\R$ is a bounded continuous function. For the equation~(\ref{eqbis}), we say that a classical solution $0\le u(t,x)\le 1$ is a transition front connecting $0$ and $1$ if there is a family of real numbers~$(X(t))_{t\in\R}$ satisfying~(\ref{gtf}).

\begin{proposition}\label{pro1bis}
Under the above assumptions, for any transition front $u$ of~\eqref{eqbis} connecting~$0$ and $1$, there holds
$$\forall\,\tau\ge0,\ \ \sup_{(t,s)\in\R^2,\,|t-s|\le\tau}|X(t)-X(s)|<+\infty.$$
\end{proposition}

\noindent{\bf{Proof.}} In the proof, $u:\R^2\to[0,1]$ denotes a transition front connecting $0$ and $1$ for~\eqref{eqbis}. It follows first from~(\ref{gtf}) and the continuity of $u$ that, for every $t\in\R$, the quantity
$$x(t)=\min\Big\{x\in\R;\ u(t,x)=\frac{1}{2}\Big\}$$
is a well defined real number such that $u(t,x(t))=1/2$. From the uniformity of the limits in~(\ref{gtf}), there is a constant $C_1\ge0$ such that
\be\label{C1}
\big|X(t)-x(t)\big|\le C_1\ \hbox{ for all }t\in\R.
\ee
On the other hand, since $u:\R^2\to[0,1]$ is bounded in $\R^2$ and $f\in L^{\infty}(\R^2\times[0,1])$, standard parabolic estimates imply that $u$ is uniformly continuous in $\R^2$. Hence, there is $\eta>0$ such that
$$\frac{1}{4}\le u(t+t',x(t))\le\frac{3}{4}\ \hbox{ for all }t\in\R\hbox{ and }t'\in[0,\eta].$$
By using again~(\ref{gtf}), one gets the existence of a constant $C_2\ge0$ such that
$$\big|X(t+t')-x(t)\big|\le C_2\ \hbox{ for all }t\in\R\hbox{ and }t'\in[0,\eta].$$
Hence,
\be\label{C2}
\big|X(t+t')-X(t)\big|\le C_1+C_2\ \hbox{ for all }t\in\R\hbox{ and }t'\in[0,\eta].
\ee
Therefore, if $\tau$ is a nonnegative real number and $k\in\N$ is such that $\tau\le k\eta$, one infers that
$$|X(t)-X(s)|\le k(C_1+C_2)$$
for all $(t,s)\in\R^2$ such that $|t-s|\le\tau$. The proof of Proposition~\ref{pro1bis} is thereby complete.~\hfill$\Box$


\end{document}